\documentclass{article}
\usepackage{amssymb,amsmath,theorem,euscript}

\input{epsf}

\newcounter{sec}

\def\sm{\smallskip}

\newcounter{punct}[sec]

\def\punct{\refstepcounter{punct}{\arabic{sec}.\arabic{punct}.  }}

\def\COUNTERS{\addtocounter{sec}{1}
              \setcounter{punct}{0}
          \setcounter{equation}{0}
          \setcounter{theorem}{0}
                  }

\newtheorem{theorem}{Theorem}[sec]
\newtheorem{proposition}[theorem]{Proposition}
\newtheorem{lemma}[theorem]{Lemma}
\newtheorem{corollary}[theorem]{Corollary}

\newtheorem{definition}[theorem]{Definition}
\newtheorem{conjecture}[theorem]{Conjecture}

\begin{document}

 \def\ov{\overline}
\def\wt{\widetilde}
 \newcommand{\rk}{\mathop {\mathrm {rk}}\nolimits}
\newcommand{\Aut}{\mathop {\mathrm {Aut}}\nolimits}
\newcommand{\Out}{\mathop {\mathrm {Out}}\nolimits}
 \newcommand{\tr}{\mathop {\mathrm {tr}}\nolimits}
  \newcommand{\diag}{\mathop {\mathrm {diag}}\nolimits}
  \newcommand{\supp}{\mathop {\mathrm {supp}}\nolimits}
  \newcommand{\indef}{\mathop {\mathrm {indef}}\nolimits}
  \newcommand{\dom}{\mathop {\mathrm {dom}}\nolimits}
  \newcommand{\im}{\mathop {\mathrm {im}}\nolimits}
 
\renewcommand{\Re}{\mathop {\mathrm {Re}}\nolimits}

\def\Br{\mathrm {Br}}

\def\SL{\mathrm {SL}}
\def\SU{\mathrm {SU}}
\def\GL{\mathrm {GL}}
\def\U{\mathrm U}
\def\OO{\mathrm O}
 \def\Sp{\mathrm {Sp}}
 \def\SO{\mathrm {SO}}
\def\SOS{\mathrm {SO}^*}
 \def\Diff{\mathrm{Diff}}
 \def\Vect{\mathfrak{Vect}}
\def\PGL{\mathrm {PGL}}
\def\PU{\mathrm {PU}}
\def\PSL{\mathrm {PSL}}
\def\Symp{\mathrm{Symp}}
\def\End{\mathrm{End}}
\def\Mor{\mathrm{Mor}}
\def\Aut{\mathrm{Aut}}
 \def\PB{\mathrm{PB}}
 \def\cA{\mathcal A}
\def\cB{\mathcal B}
\def\cC{\mathcal C}
\def\cD{\mathcal D}
\def\cE{\mathcal E}
\def\cF{\mathcal F}
\def\cG{\mathcal G}
\def\cH{\mathcal H}
\def\cJ{\mathcal J}
\def\cI{\mathcal I}
\def\cK{\mathcal K}
 \def\cL{\mathcal L}
\def\cM{\mathcal M}
\def\cN{\mathcal N}
 \def\cO{\mathcal O}
\def\cP{\mathcal P}
\def\cQ{\mathcal Q}
\def\cR{\mathcal R}
\def\cS{\mathcal S}
\def\cT{\mathcal T}
\def\cU{\mathcal U}
\def\cV{\mathcal V}
 \def\cW{\mathcal W}
\def\cX{\mathcal X}
 \def\cY{\mathcal Y}
 \def\cZ{\mathcal Z}
\def\0{{\ov 0}}
 \def\1{{\ov 1}}
 \def\frA{\mathfrak A}
 \def\frB{\mathfrak B}
\def\frC{\mathfrak C}
\def\frD{\mathfrak D}
\def\frE{\mathfrak E}
\def\frF{\mathfrak F}
\def\frG{\mathfrak G}
\def\frH{\mathfrak H}
\def\frI{\mathfrak I}
 \def\frJ{\mathfrak J}
 \def\frK{\mathfrak K}
 \def\frL{\mathfrak L}
\def\frM{\mathfrak M}
 \def\frN{\mathfrak N} \def\frO{\mathfrak O} \def\frP{\mathfrak P} \def\frQ{\mathfrak Q} \def\frR{\mathfrak R}
 \def\frS{\mathfrak S} \def\frT{\mathfrak T} \def\frU{\mathfrak U} \def\frV{\mathfrak V} \def\frW{\mathfrak W}
 \def\frX{\mathfrak X} \def\frY{\mathfrak Y} \def\frZ{\mathfrak Z} \def\fra{\mathfrak a} \def\frb{\mathfrak b}
 \def\frc{\mathfrak c} \def\frd{\mathfrak d} \def\fre{\mathfrak e} \def\frf{\mathfrak f} \def\frg{\mathfrak g}
 \def\frh{\mathfrak h} \def\fri{\mathfrak i} \def\frj{\mathfrak j} \def\frk{\mathfrak k} \def\frl{\mathfrak l}
 \def\frm{\mathfrak m} \def\frn{\mathfrak n} \def\fro{\mathfrak o} \def\frp{\mathfrak p} \def\frq{\mathfrak q}
 \def\frr{\mathfrak r} \def\frs{\mathfrak s} \def\frt{\mathfrak t} \def\fru{\mathfrak u} \def\frv{\mathfrak v}
 \def\frw{\mathfrak w} \def\frx{\mathfrak x} \def\fry{\mathfrak y} \def\frz{\mathfrak z} \def\frsp{\mathfrak{sp}}
 \def\bfa{\mathbf a} \def\bfb{\mathbf b} \def\bfc{\mathbf c} \def\bfd{\mathbf d} \def\bfe{\mathbf e} \def\bff{\mathbf f}
 \def\bfg{\mathbf g} \def\bfh{\mathbf h} \def\bfi{\mathbf i} \def\bfj{\mathbf j} \def\bfk{\mathbf k} \def\bfl{\mathbf l}
 \def\bfm{\mathbf m} \def\bfn{\mathbf n} \def\bfo{\mathbf o} \def\bfp{\mathbf p} \def\bfq{\mathbf q} \def\bfr{\mathbf r}
 \def\bfs{\mathbf s} \def\bft{\mathbf t} \def\bfu{\mathbf u} \def\bfv{\mathbf v} \def\bfw{\mathbf w} \def\bfx{\mathbf x}
 \def\bfy{\mathbf y} \def\bfz{\mathbf z} \def\bfA{\mathbf A} \def\bfB{\mathbf B} \def\bfC{\mathbf C} \def\bfD{\mathbf D}
 \def\bfE{\mathbf E} \def\bfF{\mathbf F} \def\bfG{\mathbf G} \def\bfH{\mathbf H} \def\bfI{\mathbf I} \def\bfJ{\mathbf J}
 \def\bfK{\mathbf K} \def\bfL{\mathbf L} \def\bfM{\mathbf M} \def\bfN{\mathbf N} \def\bfO{\mathbf O} \def\bfP{\mathbf P}
 \def\bfQ{\mathbf Q} \def\bfR{\mathbf R} \def\bfS{\mathbf S} \def\bfT{\mathbf T} \def\bfU{\mathbf U} \def\bfV{\mathbf V}
 \def\bfW{\mathbf W} \def\bfX{\mathbf X} \def\bfY{\mathbf Y} \def\bfZ{\mathbf Z} \def\bfw{\mathbf w}
 \def\R {{\mathbb R }} \def\C {{\mathbb C }} \def\Z{{\mathbb Z}} \def\H{{\mathbb H}} \def\K{{\mathbb K}}
 \def\N{{\mathbb N}} \def\Q{{\mathbb Q}} \def\A{{\mathbb A}} \def\T{\mathbb T} \def\P{\mathbb P} \def\G{\mathbb G}
 \def\bbA{\mathbb A} \def\bbB{\mathbb B} \def\bbD{\mathbb D} \def\bbE{\mathbb E} \def\bbF{\mathbb F} \def\bbG{\mathbb G}
 \def\bbI{\mathbb I} \def\bbJ{\mathbb J} \def\bbK{\mathbb K} \def\bbL{\mathbb L} \def\bbM{\mathbb M} \def\bbN{\mathbb N} \def\bbO{\mathbb O}
 \def\bbP{\mathbb P} \def\bbQ{\mathbb Q} \def\bbS{\mathbb S} \def\bbT{\mathbb T} \def\bbU{\mathbb U} \def\bbV{\mathbb V}
 \def\bbW{\mathbb W} \def\bbX{\mathbb X} \def\bbY{\mathbb Y} \def\kappa{\varkappa} \def\epsilon{\varepsilon}
 \def\phi{\varphi} \def\le{\leqslant} \def\ge{\geqslant}

\def\UU{\bbU}
\def\Mat{\mathrm{Mat}}
\def\tto{\rightrightarrows}

\def\Gr{\mathrm{Gr}}

\def\graph{\mathrm{graph}}

\def\O{\mathbb{O}}

\def\la{\langle}
\def\ra{\rangle}

\def\B{\mathrm B}
\def\Int{\mathrm{Int}}
\def\LGr{\mathrm{LGr}}


\def\I{\mathbb I}
\def\M{\mathbb M}
\def\T{\mathbb T}

\def\Lat{\mathrm{Lat}}
\def\LLat{\mathrm{LLat}} 
\def\Mod{\mathrm{Mod}}
\def\LMod{\mathrm{LMod}}
\def\Naz{\mathrm{Naz}}
\def\naz{\mathrm{naz}}
\def\bNaz{\mathbf{Naz}}
\def\AMod{\mathrm{AMod}}
\def\ALat{\mathrm{ALat}}

\def\Ver{\mathrm{Vert}}
\def\Bd{\mathrm{Bd}}
\def\We{\mathrm{We}}
\def\Heis{\mathrm{Heis}}

\def\bbot{{\bot\!\!\!\bot}}

\begin{center}
\Large\bf
Infinite-dimensional $p$-adic groups,
 semigroups of double cosets,
and  inner functions on Bruhat--Tits buildings

\bigskip

\large\sc
Yury A. Neretin\footnote{Supported by the grants FWF, P22122, and FWF, P25142.}
\end{center}

{\small
We construct $p$-adic analogs of operator colligations and their characteristic functions.
Consider a $p$-adic group $\bfG=\GL(\alpha+k\infty, \Q_p)$, its subgroup $L=\OO(k\infty,\Z_p)$,
and the subgroup $\bfK=\OO(\infty,\Z_p)$ embedded to $L$  diagonally. We show that double cosets
$\Gamma= \bfK\setminus \bfG/\bfK$ admit a structure of a semigroup, $\Gamma$ 
acts naturally in $\bfK$-fixed vectors
of any unitary representations of $\bfG$. For each double coset we assign a 'characteristic function',
which sends a certain Bruhat--Tits building  to another  building
(buildings are finite-dimensional); image of the distinguished  boundary
is contained in the distinguished boundary. The latter building admits a structure of (Nazarov) semigroup,
 the product in $\Gamma$ corresponds to a point-wise product of characteristic functions.
}

\section{Degeneration of Iwahori--Hecke type algebras in the infinite dimensional limit}

\COUNTERS

{\bf \punct Hypergroups of double cosets.} Consider a group $G$ and its compact subgroup $K$.
Consider double cosets $K\setminus G/K$, i.e., the quotient of $G$ with respect to the equivalence
relation
$$
g\sim k_1 g k_2,\qquad \text{where $k_1$, $k_2\in K$}
.
$$
Each double coset $\frg=KgK$ is equipped with a unique probability measure
$\mu_\frg$, which is invariant with respect to left and right translations  
by elements of $K$. Convolution of measures $\mu_{\frg}$, $\mu_{\frh}$
can be represented in the form
$$
\mu_{\frg}\ast\mu_{\frh}=\int_{K\setminus G/K} \mu_\frr(\frr) d\sigma_{\frg,\frh}(\frr)
,$$
where $\sigma_{\frg,\frh}$ is a positive probability measure on $K\setminus G/K$.
Thus we get a map
$$
(\frg,\frh)\mapsto \sigma_{\frg,\frh}
$$
from $K\setminus G/K\times K\setminus G/K$ to the space of measures
on $K\setminus G/K$. Such algebraic structures are called {\it hypergroups%
 \footnote{See, e.g., \cite{Wil}.}}.
Also the map $g\mapsto g^{-1}$ induces an involution $\mu\mapsto\mu^*$ on the hypergroup,
$$
(\mu_{\frg}\ast\mu_{\frh})^*= \mu_{\frh}^*\ast\mu_{\frg}^*
.
$$

{\sc Remark.} We reformulate this in two forms.

\sm

a)  Denote by $\cM(K\setminus G/K)$
the set of all (sign-indefinite) compactly supported measures on $G$, which are invariant with respect
to left and right translations by elements of $K$. Then $\cM(K\setminus G/K)$ is an algebra with respect
to the convolution.

b)
Let $G$ be a locally compact group with two-side invariant Haar measure
$dg$. Consider the set $C(K\setminus G/K)$ of compactly supported
left-right $K$-invariant continuous functions on $G$. Then $C(K\setminus G/K)$
is an algebra with respect to the convolution.  Sometimes it is called 
{\it (generalized) Iwahori--Hecke algebra}.
\hfill $\boxtimes$

\sm

Let $\rho$ be a unitary representation of $G$ in a Hilbert space $H$.
Denote by $H^K$ the space of $K$-fixed vectors, by $P^K$ the projection operator to $H^K$.
Let $g\in\frg$. Define the operator $H^K\to H^K$ given by
\begin{equation}
\ov\rho(\frg):=P^K \rho(g)\Bigr|_{H^K}
.
\label{eq:ov-rho}
\end{equation}
It is easy to see that $\ov \rho(\frg)$ depends on the double coset and not on a representative
$g$. 
The operators $\ov\rho(g)$ also can be expressed as
$$
\ov\rho(g)=
\int_{K\times K}\rho (k_1gk_1)dk_1\,dk_2\biggr|_{H^K}=\int_K \rho(kg)\,dk\biggr|_{H^K}
.
$$
 Also, we have a representation of the hypergroup in $H^K$ in the following sense:
$$
\ov\rho(\frg*\frh)=\int \ov\rho(\frr) d\sigma_{\frg,\frh}(\frr).
$$
Several special cases of this construction are widely used in representation theory,
in particular for the following pairs $G\supset K$:

\sm

--- $G$ is a real semisimple Lie group and $K$ is the maximal compact subgroup;
or $G$ is a compact Lie group and $K$ is a symmetric subgroup, \cite{Gel}, \cite{God}; 

\sm

--- $G$ is a finite Chevalley group, $K$ is a Borel subgroup, \cite{Iwa};

\sm

--- $G$ is a $p$-adic semisimple group and $K$ is the Iwahori subgroup, \cite{IM}.

\sm

Even
for $(G,K)=(\SL(2,\R),\SO(2))$ the explicit expression for $\sigma_{\frg,\frh}$ 
is nontrivial, see \cite{Koo}.

For smaller subgroups $K\subset G$ in semisimple groups,
the hypergroups $K\setminus G/K$ became too complicated objects. For a noncompact 
subgroup $K$ there is no {\it finite} $K\times K$-invariant measure on $K\setminus G/K$.
On the other hand, a convolution of infinite  measures is not defined
(except few exotic cases).

 In 1970s R.S.Ismagilov and G.I.Olshanski observed that  the situation can  drastically change
for  infinite-dimensional groups. Now we discuss a real archetype of our $p$-adic construction.

\sm

{\bf\punct Colligations.} Denote by $\U(\infty)$ the group
of all {\it finitary}%
\footnote{An infinite matrix $g$ is finitary, if $g-1$ has only finite number 
of nonzero matrix elements} 
infinite unitary matrices $g$ . Denote by $\OO(\infty)\subset\U(\infty)$ the group of 
real orthogonal matrices. We also use notation $\U(n+\infty)$ for the group
of block finitary unitary 
matrices $\begin{pmatrix}a&b\\c&d\end{pmatrix}$ 
of size $(n+\infty)\times(n+\infty)$. Consider double cosets 
$$
K\setminus G/K =
\OO(\infty)\setminus \U(n+\infty)/\OO(\infty),
$$
i.e., matrices $\begin{pmatrix}a&b\\c&d \end{pmatrix}\in\U(n+\infty)$ determined up to 
the equivalence
$$
\begin{pmatrix}a&b\\c&d \end{pmatrix} \sim
\begin{pmatrix}1&0\\0&u \end{pmatrix}
\begin{pmatrix}a&b\\c&d \end{pmatrix}
\begin{pmatrix}1&0\\0&v \end{pmatrix},\qquad
\text{where $u$, $v\in\OO(\infty)$}
.$$
We call such equivalence classes by {\it colligations}%
\footnote{This is a term from operator theory, a colligation (node) is
the conjugacy class (\ref{eq:conjugacy})}. 

There is no Haar measure on $K$, therefore there are no natural measures 
on double cosets $KgK$, therefore we can not repeat the construction of a hypergroup
$K\setminus G/K$.

However, there is a natural multiplication
$$
K\setminus G/K\,\times\, K\setminus G/K\,\, \to K\setminus G/K
$$
given by
\begin{equation}
\label{eq:colligations}
\begin{pmatrix}a&b\\c&d \end{pmatrix}
\circ 
\begin{pmatrix}p&q\\r&t \end{pmatrix}
=
\begin{pmatrix}a&b&0\\c&d&0\\ 0&0&1 \end{pmatrix}
\begin{pmatrix}p&0&q\\0&1&0  \\ r&0&t \end{pmatrix}
=
\begin{pmatrix}
ap&b&cq\\
cp&d&cq\\
r&0&cq
\end{pmatrix}.
\end{equation}
The matrix in the right-hand side has size $(n+\infty+\infty)$,  we regard it as a matrix of size
$$
\bigl(n+(\infty+\infty)\bigr)\times \bigl(n+(\infty+\infty)\bigr)=
(n+\infty)\times (n+\infty)
.
$$

\begin{proposition}
The $\circ$-multiplication is a well-defined associative operation on 
$K\setminus G/K$.
\end{proposition}

We also define an involution $\frg\mapsto \frg^*$ on $K\setminus G/K$
induced by the map $g\mapsto g^*$ (taking of adjoint operator).
It is easy to verify the identity
$$
(\frg\circ \frh)^*=\frh^*\circ \frg^*.
$$

Consider a unitary representation of $G=\U(n+\infty)$ in a Hilbert space
$H$. As above consider the space $H^K$ of $K$-fixed vectors in $H$ and operators
(\ref{eq:ov-rho}). The following  {\it multiplicativity theorem}
holds:

\begin{theorem}
{\rm(see \cite{Olsh-GB}, \cite{Ner-book}, Section IX.4)}
For any $\frg$, $\frh\in K\setminus G/K$,
$$
\ov\rho(\frg)\ov\rho(\frh)=\ov\rho(\frg\circ\frh).
$$
\end{theorem}

Also, for any $\frg$,
$$
\ov\rho(\frg^*)=\ov\rho(\frg)^*
$$

These phenomena (semigroup structure on $K\setminus G/K$ and the multiplicativity)
have no finite-dimensional analogs. However, for infinite-dimensional
groups they are usual, see a discussion in Subsection \ref{ss:dc}.

\sm

{\bf\punct Characteristic functions.}
We wish to describe the $\circ$-multiplication on more usual language.
For a matrix $g=\begin{pmatrix}a&b\\ c&d\end{pmatrix}$ we write the following
equation
\begin{equation}
\label{eq:char-uo}
\begin{pmatrix} 
q_+\\\lambda y\\q_-\\y
\end{pmatrix}
=
\begin{pmatrix}
\begin{pmatrix}a&b\\ c&d\end{pmatrix}&\\
&\begin{pmatrix}a&b\\ c&d\end{pmatrix}^{t-1}
\end{pmatrix}
\begin{pmatrix}
p_+\\x\\p_-\\ \lambda x
\end{pmatrix}
,
\end{equation}
where $\lambda\in\C$, $x$, $y\in\ell_2$, $p_\pm$, $q_\pm\in \C^n$.

Eliminate variables $x$, $y$ from this system of equations,
this is possible if
$$\det(\lambda^2 \ov d-d)$$
is not identical zero.
 We get
a dependence 
$$
\begin{pmatrix}
q_+\\q_-
\end{pmatrix}=
\chi_g(\lambda)
\begin{pmatrix}
p_+\\p_-
\end{pmatrix}
,$$
where $\lambda\mapsto\chi_g(\lambda)$ is a matrix-valued rational function on
$\C$. It is called a {\it characteristic function}.

A  characteristic function
 $\chi_g(\lambda)$ depends only on a double coset
$\frg$ containing $g$ and not on $g$ itself.

\begin{theorem}
If $\chi_{\frg}(\lambda)$ and $\chi_\frh(\lambda)$ are well-defined, then
\begin{equation}
\chi_{\frg\circ\frh}(\lambda)=\chi_{\frg}(\lambda) \chi_\frh(\lambda)
.
\label{eq:uo-multiplicativity}
\end{equation}
\end{theorem}

Also,
$$
\chi_{\frg^*}(\lambda)= \chi_{\frg^*}(\lambda^{-1})^{-1}.
$$


{\bf\punct Reformulation. The language of Grassmannians.}
Fix $\lambda$. Consider the set $\cX_\frg(\lambda)$ of all 
$(q_+,q_-;p_+,p_-)\in \C^{2n}\oplus \C^{2n}$ such that there are
$x$, $y$ satisfying
(\ref{eq:char-uo}).  Evidently, $\cX_\frg(\lambda)$ is a linear subspace.
Notice, that at a non-singular point of the function $\chi_\frg(\lambda)$,
the subspace
$\cX_\frg(\lambda)$ is the graph of the operator $\chi_\frg(\lambda):\C^n\to\C^n$.

Next, we  extend the function $\cX_{\frg}(\lambda)$ to the Riemann
sphere $\ov \C=\C\cup\infty$ in the following way. We write the equation 
\begin{equation}
\label{eq:char-uo-infty}
\begin{pmatrix} 
q_+\\ y\\q_-\\0
\end{pmatrix}
=
\begin{pmatrix}
\begin{pmatrix}a&b\\ c&d\end{pmatrix}&\\
&\begin{pmatrix}a&b\\ c&d\end{pmatrix}^{t-1}
\end{pmatrix}
\begin{pmatrix}
p_+\\0\\p_-\\ x
\end{pmatrix}
,
\end{equation}
and consider the set $\cX_\frg(\infty)$
 of all $(q_+,q_-;p_+,p_-)\in \C^{2n}\oplus \C^{2n}$
 such that the equation (\ref{eq:char-uo-infty})
 has a solution.

\begin{theorem}
{\rm a)} $\dim \cX_\frg(\lambda)=2n$ for all $\lambda\in\ov\C$.

\sm

{\rm b)} For any $\frg$ the map $\lambda\mapsto \cX_\frg(\lambda)$
is holomorphic on $\ov\C$.
\end{theorem}

Emphasize that the characteristic function $\cX_\frg(\lambda)$
is well-defined for all double cosets
 $\frg$.

Next, we explain how to interpret formula (\ref{eq:uo-multiplicativity})
on the language of Grassmannian. 

Let $V$, $W$ be linear spaces. We say that a {\it linear relation}
$L:V\tto W$
is a subspace $L\subset V\oplus W$.

\sm

{\sc Example.} Let $A:V\to W$ be a linear operator. Then its graph $\graph(A)\subset V\oplus W$
is a linear relation. The set of all linear subspaces in $V\oplus W$ consists
of $\dim V+\dim W$ components. Graphs of operators constitute an open dense subspace
in one of components.
\hfill $\boxtimes$

\sm

Consider two linear relations $L:V\tto W$, $M:W\tto Y$.
Define their {\it product} $LM:V\tto Y$ as the set of
$(r,p)\in V\oplus Y$ such that there exists $q\in W$
such that $(r,q)\in L$, $(q,p)\in M$. 

Also, for a linear relation $L:V\tto W$ we define
the {\it kernel} $\ker L\subset V$ and the {\it indefinity}
$\indef L\subset W$,
$$
\ker L:=L\cap(V\oplus 0),\qquad
\indef L:=L\cap(0\oplus W)
.$$

\begin{theorem}
For any $\frg$, $\frh$ and each  $\lambda\in\ov \C$,
$$
\cX_{\frg\circ\frh}(\lambda)=\cX_\frg(\lambda)\, \cX_\frh(\lambda)
.$$
\end{theorem}


{\bf \punct Conditions for characteristic functions.}
We equip the space $\C^n\oplus \C^n$ with a standard skew-symmetric
bilinear form determined by the matrix $\begin{pmatrix}0&1\\-1&0 \end{pmatrix}$.
We regard vectors $(p_+,p_-)$ and $(q_+,q_-)$ as elements of $\C^n\oplus \C^n$.
Denote by $\Sp(2n,\C)$ the group of operators preserving this form.

Equip the space $(\C^n\oplus \C^n)\oplus (\C^n\oplus \C^n)$
by the difference of skew-symmetric forms, i.e. by the form with matrix
$$
\begin{pmatrix}
 0&1&0&0\\
-1&0&0&0\\
0&0&0&-1\\
0&0&1&0
\end{pmatrix}
$$
We regard vectors $(p_+,p_-,q_+,q_-)$ as elements of this space.

\begin{proposition}
{\rm (see \cite{Ner-book}, IX.4)}

{\rm a)}
Outside poles, values of $\chi_\frg(\lambda)$
are contained in the complex symplectic group $\Sp(2n,\C)$.

\sm

{\rm b)} The characteristic function $\cX_\frg(\lambda)$ takes
values in the Lagrangian Grassmannian%
\footnote{Recall that a subspace $L$ in a $2m$-dimensional linear space
equipped with a nondegenerate skew-symmetric bilinear form is Lagrangian if
the form vanishes on $L$ and $\dim L=m$, see, e.g., \cite{Ner-gauss}, Section 3.1.}.
\end{proposition}

Second, consider the Hermitian form $M$ on $\C^n\oplus\C^n$ determined by the 
matrix $\begin{pmatrix}0&i\\-i&0 \end{pmatrix}$.
Denote by $\U(n,n)$ the group of matrices preserving $M$.

We say that a linear operator $A$ in $\C^n\oplus\C^n$
is an {\it $M$-contraction} (see, e.g., \cite{Ner-gauss}, Section 2.7), if for all vectors $v$
we have
$$
M(Av,Av)\le M(v,v)
.
$$
We say that $A$ is an {\it $M$-dilatation} if $M(Av,Av)\ge M(v,v)$.

Also, equip the space $(\C^n\oplus\C^n)\oplus(\C^n\oplus\C^n)$
with the difference of Hermitian forms, i.e. with a form
$\wt M$ given by
$$
\begin{pmatrix}
 0&i&0&0\\
-i&0&0&0\\
0&0&0&-i\\
0&0&i&0
\end{pmatrix}
$$

\begin{proposition}
\label{pr:uo-inner}
{\rm (see \cite{Ner-book}, Section IX.4)}
Let $\chi_\frg(\lambda)$ be well-defined. Then:

\sm

{\rm a)} If $|\lambda|=1$, then $\chi_\frg(\lambda)\in\U(n,n)$.

\sm
 
{\rm b)} If $|\lambda|<1$, then $\chi_\frg(\lambda)$ is 
an $M$-contraction.

\sm

{\rm c)} If $|\lambda|>1$, then $\chi_\frg(\lambda)$ 
is an $M$-dilatation.
\end{proposition}

\begin{proposition}
\label{pr:uo-inner-1}
{\rm (see \cite{Ner-book}, Section IX.4)}

{\rm a)} If $|\lambda|=1$, then the subspace $\cX_\frg(\lambda)$
is $\wt M$-isotropic.

\sm
 
{\rm b)} If $|\lambda|<1$, then the form
$\wt M$ is positive semi-definite on the subspace $\cX_\frg(\lambda)$.

\sm

{\rm c)}  If $|\lambda|>1$,
then the form
$\wt M$ is negative semi-definite on the subspace $\cX_\frg(\lambda)$.

\sm

{\rm d)} If $|\lambda|<1$, then the from $M$ is strictly positive definite on%
\footnote{This condition contains additional information  
only at points $\lambda$,  where $\cX(\lambda)$ is not a graph 
of an operator. By statement b) $M$ is positive semi-definite on the kernel.}
 $\ker\cX_\frg(\lambda)$. Also $M$ is negative definite on $\indef \cX_\frg(\lambda)$.
\end{proposition}

Characteristic functions also satisfy to the following condition of symmetry
at 0
\begin{equation}
\chi_\frg(-\lambda)=
\begin{pmatrix}1&0\\0&-1\end{pmatrix} \chi_\frg(\lambda)
 \begin{pmatrix}1&0\\0&-1\end{pmatrix}^{-1}
.
\label{eq:char-uo-last}
\end{equation}
On the language of Grassmannians this means
\begin{equation}
(p_+,p_-,q_+,q_-)\in\cX(\lambda)\Leftrightarrow (p_+,-p_-,q_+,-q_-)\in\cX(\lambda)
.
\label{eq:char-uo-last-1}
\end{equation}

\begin{theorem}
Any holomorphic map  $\cX$ from $\ov\C$ to the
Lagrangian Grassmannian   satisfying
the conditions of
Proposition {\rm\ref{pr:uo-inner-1}} and condition 
{\rm(\ref{eq:char-uo-last-1})} is a
characteristic function of a double coset $\frg$.
\end{theorem}


{\bf\punct Central extension.} 
A characteristic function is not sufficient for a
reconstruction of a double coset, in fact matrices of the form
$$
\left(
\begin{matrix}a&b&0&\quad \}n \\ c&d&0&\quad \}\infty \\0&0&e&\quad \}\infty \end{matrix}
\!\!\!\!\!\!\!\!\!\!\!\!\!\!\!\!\!\!
\right)
$$
with fixed
 $a$, $b$, $c$, $d$ and arbitrary  
  $e$  have the same characteristic function.
Let us introduce an additional invariant. We write
the equation
 $$
\begin{pmatrix} 
0\\\lambda y\\0\\y
\end{pmatrix}
=
\begin{pmatrix}
\begin{pmatrix}a&b\\ c&d\end{pmatrix}&\\
&\begin{pmatrix}a&b\\ c&d\end{pmatrix}^{t-1}
\end{pmatrix}
\begin{pmatrix}
0\\x\\0\\ \lambda x
\end{pmatrix}
,$$
as an equation for $x$, $y$. Denote by $n_\frg(\lambda)$
the dimension of the space of solutions of this equation.
Then 

\sm

--- $n_\frg(\lambda)=0$ for all but a finite number of values of $\lambda$;

\sm

--- $n_\frg(\lambda)=0$ if $|\lambda|\ne 1$;

\sm

--- $n_\frg(\lambda)=n_\frg(-\lambda)$;

\sm

--- $n_\frg(\pm 1)=\infty$.

\sm

Thus we get a finite set with multiplicities (we call it {\it divisor}).

\begin{theorem}%
\footnote{This and previous statements are
given in  \cite{Ner-book},
я.IX.4.8 without formal proof. In fact, a
proof is contained in the same book, Addendum E.
Precisely, in Subsection
E.4 it is shown how to reduce our statements to the standard theorem
 (see \cite{Bro}) 'pure unitary operator node is determined by
 its characteristic function'. In fact, we only need   this theorem 
 for finitary matrices and rational characteristic functions.} 
A double coset is uniquely determined by its characteristic function $\cX$
and the divisor $n$.  
\end{theorem}

\begin{theorem} {\rm \cite{Ner-book}, IX.4.5)}
$$
n_{\frg\circ\frh}(\lambda)=
n_{\frg}(\lambda)+ n_\frg(\frh;\lambda)+
\dim \Bigl(\indef \cX_\frh(\lambda)\cap \ker \cX_\frg(\lambda)\Bigr)
.
$$
\end{theorem}
 
Double cosets corresponding matrices 
$$
\left(
\begin{matrix}1&0&0&\quad \}n \\ 0&1&0&\quad \{\infty \\0&0&e&\quad \}\infty \end{matrix}
\!\!\!\!\!\!\!\!\!\!\!\!\!\!\!\!\!\!
\right)
$$
is the center of the semigroup $K\setminus G/K$. The quotient
of $K\setminus G/K$ with respect to the center is isomorphic to the semigroup
of rational matrix-values functions described above. 

\sm

{\bf\punct Degeneration of hypergroups of double cosets.}
Let $N>k$. Embed $\U(n+k)$ to $\U(n+k+N)$ by
 $$
 \iota_N:\begin{pmatrix}A&B\\C&D\end{pmatrix}
 \mapsto 
 \begin{pmatrix}A&B&0\\C&D&0\\0&0&1\end{pmatrix}
 .
 $$
 Embed $\U(k+N)$ to $\U(n+k+N)$ by
$$
\begin{pmatrix}\alpha&\beta\\ \gamma&\delta \end{pmatrix}
\mapsto \begin{pmatrix}1&0&0\\0& \alpha&\beta\\ 0&\gamma&\delta \end{pmatrix}
.
$$
Fix matrices $g=\begin{pmatrix}q&b\\c&d \end{pmatrix}$,
 $h=\begin{pmatrix} p&q\\r&t \end{pmatrix}\in \U(n+k)$. 
 Then  for $N>k$ a matrix
 $\iota_N(g)\circ \iota_N(h)$ is well-defined as an element
 of $\U(k+N) \setminus \U(n+k+N)/\U(k+N)$.

We equip the group $\U(n+k+N)$ with the metric induced by the operator norm
in Euclidean $\C^{n+k+N}$.

\begin{proposition}
\label{pr:concentration}
 Fix $g$, $h\in \U(n+k)$ as above.
Consider the corresponding double cosets 
$$\frg_N, \frh_N\in \U(k+N) \setminus \U(n+k+N)/\U(k+N)$$
and the measure 
$$
\kappa_N=\mu_{\frg_N}*\mu_{\frh_N}
$$ 
Then for each $\epsilon>0$, $\delta>0$ there exists $N$ such that
the measure $\kappa_N$ of $\epsilon$-neighborhood of 
$\iota_N(g)\circ \iota_N(h)$ is $>1-\delta$.
\end{proposition}

See \cite{Olsh-GB}, \cite{Olsh-topics}, \cite{Ner-haar}.

\sm
 
{\bf\punct Semigroups of double cosets.%
\label{ss:dc}}
The first example of multiplication of double cosets
 was discovered by Ismagilov \cite{Ism1},
 he considered the group $G=\SL(2,k)$ over a non-Archimedian normed non locally compact field
 $k$.
The subgroup $K$ is the group $\SL(2,o)$ over integer elements of $k$.
The double cosets are parametrized by non-negative integers $\Z_+$, and the operation 
$\circ$
is the usual addition. The multiplicativity theorem allows  to classify spherical functions
(see also \cite{Ism2}).
 Olshanski
\cite{Olsh-tree}  showed that this semigroup
is  a limit of hypergroups $\SL(2,\Z_p)\setminus\SL(2,\Q_p)/ \SL(2,\Z_p)$
as $p\to\infty$.

Next, consider a series of Riemannian symmetric spaces $G(n)/K(n)$
(an example is $\U(n)/\OO(n)$). Olshanski \cite{Olsh-GB}, \cite{Olsh-topics} showed that
the same phenomena hold for any pair $G(k+\infty)\supset K(\infty)$. Also he described such semigroups
for infinite symmetric groups. As far as we know description of such objects, they became
a tool of the representation theory. On the other hand, it seems that such structure are
interesting by themselves.

In \cite{Ner-book}, Section 8.5, the author  observed
that  multiplications on $K\setminus G/K$ are quite usual for infinite-dimensional
groups (see also \cite{Ner-faa}, \cite{Ner-fein}).
  In fact this happened more-or-less always if $K$ is one of the following
groups:

\sm

1) $K$ is a complete infinite unitary group, orthogonal group, or
symplectic (quaterninic unitary) group (or a product of several copies
of such groups);

\sm

2) $K$ is the infinite symmetric group $S(\infty)$;

\sm

3) $K$ is the group of automorphisms of a measure space;

\sm

These groups are infinite-dimensional imitation of compact groups (but they are neither compact,
nor locally compact)
apparently some other examples also exist (for instance, below we discuss
$K=\OO(\infty,\Z_p)$).

For precise general theorems, see \cite{Ner-faa}, \cite{Ner-fein}.
To explore them  we need explicit descriptions of $K\setminus G/K$,
  such descriptions  recently were obtained in \cite{Ner-symm}, \cite{Ner-char}, 
  \cite{Ner-faa}, \cite{Ner-fein}.







\sm

{\bf\punct Inner functions.} Recall a definition of inner functions.

\sm

1) A holomorphic function $f(z)$ in a unit disk $|z|<1$
is called {\it inner},  if $|f(z)|<1$ for $|z|<1$ and
\begin{equation}
\lim_{r\to 1^-} |f(r e^{i\theta})|=1
\qquad\text{a.s. $\theta\in [0,2\pi]$}
\label{eq:radial}
,
\end{equation}
where $z=r e^{i\theta}$ and $r$, $\theta$ are real%
\footnote{We can not write a limit as  $z\to e^{i\theta}$, an inner function can be discontinuous
at all points of the circle.}. 
On this topic, see, e.g., \cite{Garn}. It can be shown
that limit (\ref{eq:radial}) can be replaced a.s. by the
{\it nontangential limit}
\begin{equation}
\lim\limits_{z\to e^{i\theta},\,\,
\left|\arg\frac{e^{i\theta}-z}{e^{i\theta}}\right|\le \pi/2-\epsilon} f(z)
,\end{equation}
where $\epsilon>0$ is fixed (in fact we consider a limit over the angle
whose vertex is $e^{i\theta}$, the bisector is $t e^{i\theta}$, and
the value of the angle is $\pi-2\epsilon$.

 \sm
 
2) A homomorphic  matrix-valued (operator-valued) function 
 $f(z)$ in the unit disk
is called {\it inner} if $\|f(z)\|\le 1$ for $|z|<1$ and boundary values    
of $f$ on the circle are unitary  (see Livshits \cite{Liv1}, Potapov \cite{Pot}).
 Consider an  operator $d$  closed to unitary 
(one of possible variants $\rk(dd^*-1)=\rk (d^*d-1)<\infty$) with $\|d\|=1$.
We are interested its properties up to conjugations $d\mapsto udu^{-1}$, 
where $u$ is unitary.
Build a larger {\it unitary} matrix
$g=\begin{pmatrix}a&b\\c&d \end{pmatrix}$ including $d$ as a block. We consider
$g$ up to the equivalence 
\begin{equation}
\begin{pmatrix}a&b\\c&d \end{pmatrix} \sim
\begin{pmatrix}1&0\\0&u \end{pmatrix}
\begin{pmatrix}a&b\\c&d \end{pmatrix}
\begin{pmatrix}1&0\\0&u^{-1} \end{pmatrix},\qquad
\text{where $u\in\U(\infty)$}
.
\label{eq:conjugacy}
\end{equation}
Assign to $g$ the expression ({\it characteristic function}) by
$$
\chi(\lambda)=
a+\lambda b(1-\lambda d)^{-1}c
.
$$
 Such functions (under some conditions on $d$) are inner functions $\theta(z)$ in the unit disk.
 Invariant subspaces of $d$
are in one-to-one correspondence with divisors of $\theta$ in the class of
inner functions. The product of inner functions 
corresponds to the product of conjugacy classes (\ref{eq:conjugacy})
by formula (\ref{eq:colligations}).

\sm

3) More generally, consider a pseudo-Euclidean space. 
We say that a meromorphic matrix-valued function $f$ in the disk is {\it inner} if it is
indefinite contractive in the disk and pseudo-unitary on the unit circle.
Such functions arise in the same context but the condition $\|d \|\le 1$
is omitted.

The characteristic function of double cosets defined  above are inner in this sense.

\sm

4) Denote by $B_n$ the set of all $n\times n$ complex symmetric matrix
with norm $< 1$; $B_n$ also is an Hermitian symmetric space 
$$B_n=\U(n,n)/\U(n)\times \U(n),$$
its distinguished boundary (Shilov boundary) consist of unitary matrices.

 In   \cite{Ner-char}, \cite{Ner-faa} there were considered various 
 semigroups of double cosets on infinite-dimensional classical groups.
 For instance, consider group $G=\U(\alpha+k\infty)$ consisting of block
 unitary matrices of size $\alpha+\infty+\dots+\infty$. Consider its subgroup 
 $L=\U(\infty)$ embedded to $G$ in the block diagonal way. Consider 
 the subgroup $K=\OO(\infty)\subset L$ embedded to $\U(\infty)$ in the 
 natural way. Then $K\setminus G/K$ is a semigroup. 
 Characteristic functions \cite{Ner-char} are inner functions in $B_k\times B_k$
 taking values at the space of $2\alpha\times 2\alpha$-matrices.
 This means that values of a function are $M$-contractions inside 
 $B_k\times B_k$ and are pseudounitary on the Shilov boundary $\U(n)\times\U(n)$.
 The product of double cosets corresponds to the product of characteristic functions.
 
 It is possible to vary the definition and to regard a characteristic function
 as a map $B_{k}\times B_k\to B_{2\alpha}$.

\sm


{\bf\punct Infinite-dimensional  $p$-adic groups.}
Representation theory of
infinite-dimensional classical groups (see, e.g., \cite{Olsh-kir},
\cite{Olsh-short}, \cite{Olsh-GB}, \cite{Ner-book}, \cite{Olsh-add2},
\cite{BO}, \cite{Ner-faa})
 and infinite symmetric groups
 (see, e.g., \cite{Olsh-symm}, \cite{KOV}, \cite{Ner-symm}) exists and is well-developed.
  There were several recent works concerning infinite-dimensional
 classical groups over finite fields (see \cite{KV}, \cite{Ner1}, \cite{GKV}).

 Few is known about infinite-dimensional  $p$-adic groups. There are the following works:
 
 \sm
 
 1) Work of Nazarov \cite{Naz}, \cite{NNO} on the
 Weil representation of an infinite-dimensional group $\Sp(2\infty,\Q_p)$.
  Existence of such
 representation is more-or-less evident. However, the Weil representation
 of $\Sp(2n,\R)$ and $\Sp(2\infty,\R)$ admits a continuation
 to a certain complex domain $\Gamma$ (if $n<\infty$, then $\Gamma$ 
 is a semigroup parametrized by
 complex symmetric $2n\times 2n$ matrices with norm $<1$, 
 see, e.g., \cite{Ner-book}, Section 4.2,
 \cite{Ner-gauss}, Section 5.1).
 Nazarov constructed an analog of $\Gamma$
 for $p$-adic case, see below Section 3 (for more details, see \cite{Ner-gauss}, Sections 10.7, 11.2)

\sm

2) A construction of Hua measures on $p$-adic Grassmannians
and on the inverse limit of $p$-adic Grassmannians
in \cite{Ner-huap}. This is an analog of inverse limits of
compact symmetric spaces (see \cite{Ner-hua}) and of symmetric groups (see \cite{KOV}).
Recall that in latter two cases there exists a substantial harmonic analysis on such 
inverse limits, see \cite{BO}, \cite{KOV}.

\sm

3) The group of diffeomorphisms of $p$-adic projective line is an object similar
to the group of diffeomorphisms of the circle (many constructions of representations of
the latter group survive
in $p$-adic case, \cite{Ner-trees}).

\sm

{\bf\punct A $p$-adic example.} Here we briefly discuss a $p$-adic object,
which is related to the topic of this paper but more simple. Let $\Q_p$ be a $p$-adic
field, $\Z_p\subset \Q_p$ be the ring of $p$-adic integers. Denote by
$\GL(\infty,\Q_p)$ the group of finitary invertible matrices over $\Q_p$.
Consider conjugacy classes of $\GL(\alpha+\infty,\Q_p)$ with respect to the subgroup
$\GL(\infty,\Z_p)$,
\begin{equation}
\begin{pmatrix}a&b\\c&d \end{pmatrix} \sim
\begin{pmatrix}1&0\\0&u \end{pmatrix}
\begin{pmatrix}a&b\\c&d \end{pmatrix}
\begin{pmatrix}1&0\\0&u^{-1} \end{pmatrix},\qquad
\text{where $u\in\GL(\infty,\Z_p)$}
.
\label{eq:conjugacy-p}
\end{equation}
Such conjugacy classes admit a natural $\circ$-multiplication by formula
(\ref{eq:colligations}), this multiplication is a well-defined associative
operation on the space of conjugacy classes. We wish to construct
an analog of characteristic functions. 

First, choose a sufficiently large $m$ such that a matrix 
$\begin{pmatrix}a&b\\c&d \end{pmatrix}$ is actually contained
in $\GL(\alpha+m,\Q_p)$.
Consider a lattice%
\footnote{see a definition in Subsection \ref{ss:modules}} $R\subset \Q_p^2$.
For this lattice we consider the lattice
$$
R\otimes \Z_p^m\subset \Q_p^2\otimes \Q_p^m\simeq \Q_p^m\oplus \Q_p^m
.$$ 
We write an equation
\begin{equation}
\begin{pmatrix}
v\\y
\end{pmatrix}=
\begin{pmatrix}a&b\\c&d \end{pmatrix}
\begin{pmatrix}
u\\x
\end{pmatrix}
\label{eq:for-p}
.\end{equation}
Next, consider the set $\chi(R)$ of all pairs $(v,u)\in \Q_p^\alpha \oplus \Q_p^\alpha$
for which there exists $y\oplus x\in R\otimes \Z_p^m$
such that the equality (\ref{eq:for-p}) is satisfied.
Then $\chi(R)$ is a $\Z_p$-submodule in $\Q_p^\alpha \oplus \Q_p^\alpha$,
which can be regarded as a relation   $\Q_p^\alpha \tto \Q_p^\alpha$.
The $\circ$-product corresponds to point-wise product of functions 
$\chi(R)$ with values in relations.

We also point out that these functions are compatible with the structure of
Bruhat--Tits buildings and are inner in a reasonable sense. Both phenomena
are discussed below for more sophisticated objects.

\sm

{\bf\punct Purpose of the paper.} We wish to describe multiplication  of double cosets
on $p$-adic groups and to obtain analogs of characteristic functions.
For a double coset we assign a simplicial map from a Bruhat--Tits building
$\Omega$
to a Bruhat--Tits building $\Xi$ such that the image of the distinguished boundary is 
contained in the distinguished boundary. We also have a structure of a semigroup on the set
of vertices of the building
$\Xi$ (the Nazarov semigroup) and the product of double cosets 
corresponds to pointwise product of functions $\Omega\to \Xi$. 

 Our construction
is not a final solution of the problem%
\footnote{
 First, we do not introduce an analog of the 'divisor'. Secondly,
\cite{Ner-invariant} suggests that complete data separating double cosets must contain a sequence of characteristic functions determined on the increasing sequence of buildings. A construction of complete data in a real case in \cite{Ner-invariant} is based on classical invariant theory, which does not valid over the ring $\Z_p$.}

\sm

{\bf \punct A non-properly understood link.}
In fact our main construction below is organized as an extension of rational maps
of $p$-adic Grassmannians to simplicial maps of Bruhat--Tits buildings.
Also, our construction admits an automatic pass to algebraic extensions.
 Constructions of such
type are investigated in theory of Berkovich analytic spaces, see, e.g.,
\cite{Bak}, \cite{Con}. However their extensions are rigid, and our extensions depend
on additional data%
\footnote{Below rational maps of Grassmannians originate from double
cosets
$$\OO(\infty,\Q_p)\setminus\GL(\alpha+k\infty,\Q_p)/\OO(\infty,\Q_p)$$
(see  Proposition \ref{pr:another-quotient})
maps of Bruhat--Tits buildings from double cosets
$$\OO(\infty,\Z_p)\setminus\GL(\alpha+k\infty,\Q_p)/\OO(\infty,\Z_p).$$
Therefore we get many maps of Bruhat--Tits buildings with the same
restriction to a distinguished boundary, i.e., to the Grassmannian.}.
 So I can not understand relations of our constructions and Berkovich theory.

\sm 


{\bf \punct Notation.} Let

\sm

--- $A^t$ be the transposed matrix;

\sm

--- $1_\alpha$, $1_V$ be the unit matrix of order $\alpha$, the unit operator in a space
$V$;

\sm

--- $\Q_p$ be the $p$-adic field;

\sm

--- $\Z_p$ be the ring of 
$p$-adic integers;

\sm

--- $\Q_p^\times$, $\C^\times$ be multiplicative groups of $\Q_p$, $\C$.

\sm

We denote the standard character $\Q_p\to\C^\times$ by $\exp\{2\pi i a\}$.
For $a=\sum_{\ge-N} a_j p^j$, where $a_j=0$, $1$, \dots, $p-1$, we set
$$
\exp\{2\pi i a\}=\exp\Bigl\{2\pi i\sum_{j\ge -N} a_jp^j\Bigr\}:=
\exp\Bigl\{2\pi i\sum_{j: -1\ge j\ge -N} a_jp^j\Bigr\}
$$
Below we define:

\sm

--- the groups $\GL(n,\Q_p)$,  $\Sp(2n,\Q_p)$, $\Sp(2n,\Q_p)$,  
$\OO(n,\Z_p)$, $\GL(\infty,\Z_p)$, $\Sp(2\infty,\Q_p)$, etc.,  Subsection \ref{ss:groups};

\sm

--- $V_\pm$, formula (\ref{eq:V-pm});

\sm

--- groups $\bfG=\GL(\alpha+k\infty,\Q_p)$, $\bfK=\OO(\infty,\Z_p)$,  Subsection \ref{ss:double-cosets};

\sm

--- $\frg\star \frh$, the product of double cosets, Subsection \ref{ss:double-cosets};

\sm

--- $\frg^*$, the involution on double cosets, Subsection \ref{ss:involution};

\sm

--- 
    $R_\downarrow$, $R^\uparrow$, , Subsection \ref{ss:modules};

\sm

--- $R_j\nearrow R$, rigid convergence, \ref{ss:convergence-modules};

\sm

--- $\LMod(V)$, $\LLat(V)$, $\LGr(V)$, spaces of Lagrangian submodules, Subsection \ref{ss:self-dual};

\sm

--- $\Delta(V)$, $\Bd(V)$, buildings, Subsections \ref{ss:almost}, \ref{ss:buildings};

\sm

--- $P:V\tto W$, $\ker P$, $\indef P$, $\dom P$, $\im P$,  $P^\square$, Subsection \ref{ss:relations};

\sm

--- $\Naz$, $\ov\Naz$, $\bNaz$,  the Nazarov category, Subsections \ref{ss:nazarov};
\ref{ss:extended-nazarov};

\sm

--- $\We$, the Weil representation, Subsection \ref{ss:weil};

\sm

--- $\chi_\frg(Q,T)$, a characteristic function, Subsection 4.1.


\section{Multiplication of double cosets}

\COUNTERS

{\bf \punct Groups.%
\label{ss:groups}} By  $V=\Q_p^n$ we denote linear spaces over $\Q_p$.
Denote by $\GL(n,\Q_p)=\GL(V)$ the group of invertible linear operators in $\Q_p^n$;
by  $\GL(n,\Z_p)$ the group of all matrices $g$ with integer elements,
such that $g^{-1}$ have integer elements.

Consider a space $V=\Q_p^{2n}$ equipped with a non-degenerate skew-symmetric bilinear form
$B_V$,
say $\begin{pmatrix} 0&1\\-1&0  \end{pmatrix}$.
The symplectic group $\Sp(2n,\Q_p)$ is the group of matrices preserving this form,
$\Sp(2n,\Z_p)$ is the group of symplectic matrices with integer elements.
We also denote
\begin{equation}
V_+:=\Q_p^n\oplus 0,\qquad V_-=0\oplus \Q_p^n
.
\label{eq:V-pm}
\end{equation}

Also, consider a space $\Q_p^n$ equipped with the standard symmetric
bilinear form 
$
(v,w)=\sum v_j w_j
.$ 
We denote by $\OO(n, \Q_p)$ the group of all matrices preserving this form%
\footnote{There are several  non-equivalent non-degenerate quadratic forms
and several different orthogonal groups  over 
$\Q_p$, however we consider only this group.}.

By $\GL(\infty,\Q_p)$ we denote the group of all infinite invertible matrices
over $\Q_p$ such that $g-1$ has only finite number of non-zero elements.
We call such matrices {\it finitary}. 
We define $\GL(\infty,\Z_p)$, $\Sp(2\infty,\Q_p)$, $\Sp(2\infty,\Z_p)$,
$\OO(\infty,\Z_p)$
in the same way.

\sm


{\bf\punct Multiplication of double cosets.%
\label{ss:double-cosets}} Let 
$$\bfG:=\GL(\infty,\Q_p):=\GL(\alpha+k\infty,\Q_p)$$
 be the group of 
finitary block 
$(\alpha+\infty+\dots+\infty) \times (\alpha+\infty+\dots+\infty)$-
matrices (there are $k$ copies of $\infty$). 
By $\bfK$ we denote the group 
$$\bfK=\OO(\infty,\Z_p)$$
 embedded to
$\bfG$ by the rule
\begin{equation}
\frI: u\mapsto \begin{pmatrix}
                1_\alpha &0&\dots&o\\
                0&u&\dots &0\\
                \vdots&\vdots&\ddots&\vdots\\
                0&0&\dots &u
               \end{pmatrix}
               \label{eq:tau}
,\end{equation}
where $1_\alpha$ denotes the unit matrix of order $\alpha$.

\sm

{\sc Remark.} Certainly, $\bfG:=\GL(\infty,\Q_p)$. But the notation of the type
$\bfG:=\GL(\alpha+k\infty,\Q_p)$ allows us to indicate certain subgroups in $\bfG$.
\hfill $\boxtimes$

\sm

 We wish to define a structure of a semigroup 
on double cosets 
$\bfK\setminus \bfG/\bfK$.

Set
\begin{equation}
\Theta_N:=
\begin{pmatrix}
0&1_N&0\\
1_N&0&0\\
0&0&1_\infty
\end{pmatrix} 
\in\bfK
.
\label{eq:Theta}
\end{equation}
Let $\frg$, $\frh\in \bfK\setminus \bfG/\bfK$.
Choose their representatives $g$, $h\in \bfG$.
Consider the sequence
$$
f_N:=g\frI(\Theta_N)h
$$
and double coset $\frf_N$ containing $f_N$.

\begin{theorem}
\label{th:umnozhenie}
{\rm a)} The sequence $\frf_N$ of double cosets  is eventually constant.

\sm

{\rm b)} The limit $\frf:=\lim_{N\to\infty}\frf_N$
does not depend on a choice of representatives $g$, $h$. 

\sm 

{\bf c)} The product $\frg\star\frh$ in $\bfK\setminus \bfG/\bfK$ obtained in this way is associative. 

\end{theorem}

These statements are simple, see proofs of parallel real statements
in \cite{Ner-faa}. Also, it is easy to write an explicit formula for the product.
For definiteness, set 
$k=2$. Then 
\begin{multline*}
\begin{pmatrix} 
 a& b_1&b_2\\
c_1&d_{11}&d_{12}\\
c_2&d_{21}&d_{22}
\end{pmatrix}
\star
\begin{pmatrix} 
 a'& b'_1&b'_2\\
c'_1&d'_{11}&d'_{12}\\
c'_2&d'_{21}&d'_{22}
\end{pmatrix}
=\\=
\begin{pmatrix} 
 a& b_1&0&b_2&0\\
c_1&d_{11}&0&d_{12}&0\\
0&0&1&0&0\\
c_2&d_{21}&0&d_{22}&0\\
0&0&0&0&1
\end{pmatrix}
\begin{pmatrix}
1_\alpha&0&0&0&0\\ 
0&0&1_\infty&0&0\\
0&1_\infty&0&0&0\\
0&0&0&0&1_\infty\\
0&0&0&1_\infty&0\\
\end{pmatrix}
\begin{pmatrix} 
 a'& b'_1&0&b'_2&0\\
c'_1&d'_{11}&0&d'_{12}&0\\
0&0&1&0&0\\
c'_2&d'_{21}&0&d'_{22}&0\\
0&0&0&0&1
\end{pmatrix}
\end{multline*}
Since a result is double coset, 
we can write the final matrix in different forms,
say
\begin{equation}
\frf=
\begin{pmatrix}
aa'& |& b_1& ab_1'& & b_1& ab_1'
\\
-& + & -& -&-& - &-
\\
c_1a'&|&  d_{11}& c_1b_1'&& d_{12}& c_1b_2'
\\
c'_1&|&   0 &  d_{11}'&& 0 &  d_{12}'& 
\\
&|&&&&
\\
c_2a'&|&  d_{21}& c_2b_1'&& d_{22}& c_2b_2'
\\
c'_2& |&  0 &  d_{21}'&& 0 &  d_{22}'& 
\end{pmatrix}
\label{eq:coset-product}
\end{equation}
or
$$
\frf=
\begin{pmatrix}
aa'& |&  ab_1' & b_1&&  ab_2' & b_2
\\
-& + & -& -& -&-
\\
c_1 a'&|&  c_1 b_1'& d_{11} && c_1 b_2'& d_{12} 
\\
c_1' &|&  d_{11}'&0 &&  d_{12}'&0 
\\
&|&&&&
\\
c_2 a'&|&  c_2 b_1'& d_{21} && c_2 b_2'& d_{22} 
\\
c_2' &|&  d_{21}'&0 &&  d_{22}'&0 
\end{pmatrix}
.
$$


{\bf\punct Multiplicativity theorem.%
\label{ss:multiplicativity}}
Let $\rho$ be a unitary representation of $\bfG$, denote by $H^\bfK$ 
the subspace of all $\bfK$-fixed vectors. Denote by
$P^\bfK$ the operator of orthogonal projection
to $H^\bfK$. For $g\in \bfG$ consider the operator
$\ov\rho(g):H^\bfK\to H^\bfK$ given by
$$
\ov\rho(g):=P^\bfK \rho(g)\Bigr|_{H^K}.
$$
Obviously,
 $\ov\rho(g)$ is a function on double cosets $\bfK\setminus \bfG/\bfK$,
therefore we can write $\ov\rho(\frg)$. 

\begin{theorem}
\label{th:multiplicativity} For any unitary representation $\rho$,
for all $\frg$, $\frh\in \bfK\setminus \bfG/\bfK$ the following equality
{\rm(}the ``multiplicativity theorem''{\rm)}
holds,
$$
\ov\rho(\frg)\ov\rho(\frh)=
\ov\rho(\frg\star \frh)
.
$$
\end{theorem}

We give a proof in Section 6.

\sm

{\sc Remark.} Apparently the analog of Proposition
\ref{pr:concentration} for $p$-adic case does not hold.
\hfill $\square$

\sm


{\bf\punct Sphericity.%
\label{ss:sphericity}}

\begin{proposition}
\label{pr:spherical}
Let $\alpha=0$.
Then the pair $(\bfG,\bfK)$ is spherical, i.e., for any irreducible unitary
representation of $\bfG$ the dimension of the space of $\bfK$-fixed vectors 
is $\le 1$.
\end{proposition}

We omit a proof, it is the same as for infinite-dimensional real classical groups,
see \cite{Ner-faa}.

\sm


{\bf\punct Involution.%
\label{ss:involution}} The map $g\mapsto g^{-1}$ induces an involution
$\frg\mapsto \frg^*$ on $\bfK\setminus \bfG/\bfK$. Evidently,
$$
(\frg\star\frh)^*=\frh^*\star \frg^*
.$$
Also, for any unitary representation $\rho$ of $\bfG$ we have
$$
\ov\rho(\frg^*)=\ov\rho(\frg)^*
.
$$


\sm

{\bf\punct Purpose of the work.}
Our aim is to describe this multiplication in 
more usual terms. More precisely, we wish to get $p$-adic analogs
of multivariate characteristic functions constructed in
 \cite{Ner-faa}, \cite{Ner-char}.


\sm

{\bf \punct Structure of the paper.}
 Section 3
contains preliminaries (lattices, Bruhat--Tits buildings, relations,
the Weil representation of the Nazarov category). A main construction 
(characteristic functions of  double cosets and their properties) is contained in Section 4.
Proofs are given in Section 5.

 In  Section 6 we prove the multiplicativity
theorem. Section 7 contains some  constructions of 
representations. Theorem \ref{th:link} shows a link between the characteristic function
and operators $\ov\rho(\frg)$.

\section{Preliminaries. Submodules, relations, Bruhat--Tits  buildings,  Nazarov category,
and Weil representation}

\COUNTERS

{\bf A. Submodules and convergence}

\bigskip

{\bf \punct Modules.%
\label{ss:modules}}
Below the term {\it submodule} means an
$\Z_p$-submodule  in a linear space $V=\Q_p^k$. For each submodule $R\subset \Q_p^k$
 there is a (non-canonical) basis
$e_i\in \Q_p^k$ such that 
\begin{equation}
R=\Q_p e_1\oplus\dots \oplus \Q_p e_j\oplus  
\Z_p e_{j+1}\oplus\dots \oplus \Z_p e_l
.
\label{eq:R}
\end{equation}
If $j=k$ then $R$ is a linear subspace. If
$j=0$, $l=k$, then we get a {\it lattice}.
A formal definition is:
a {\it lattice } $R$ is a compact $\Z_p$-submodule
such that $\Q_p R=\Q_p^k$. For details,
see, e.g., \cite{Weil}.

Denote by $\Mod(V)$ the set of all submodules in $V$, by $\Lat(V)$
the space of all lattices. It is easy to see that
$$
\Lat(V)\simeq \GL(V,\Q_p)/\GL(V,\Z_p)
.$$

For any  submodule $R$ denote by $R_\downarrow$
the maximal linear subspace in $R$. By $R^\uparrow$ we denote
the minimal linear subspace containing $R$,
$$
R_\downarrow \subset R\subset R^\uparrow
$$
The image of $R$ in the quotient space $R^\uparrow/R_\downarrow$ is a lattice.

Conversely, let $L\subset M$ be a pair of subspaces,  $\pi:L\to L/M$
be the projection.  Let
$P\subset M/L$ be a lattice.
 Then $\pi^{-1} P$ is a submodule in $\Q_p^k$
and all submodules have such form.


\sm

{\bf \punct Duality.} For a $p$-adic linear space $V$ we denote by
$V'$ the space of linear functionals on $V$. For a submodule
$L\subset V$ define the dual module $L^\lozenge\subset V'$ as the set of all linear functionals
$\ell\in V'$ such that
$$
\ell(v)\in\Z_p\qquad \text{for all $ v\in L$}
$$
Notice that $L^{\lozenge\lozenge}=L$. 

If $L$ is a lattice, then $L^\lozenge$ is a lattice.

\sm

{\bf \punct The Hausdorff convergence on $\Mod(V)$.%
} Let $V=\Q_p^n$.
We define a norm on $V$ as 
$$
\|x\|=\max_{j} |x_j|
.$$
Denote by $B(p^l)$ the ball with center at 0 of radius $p^l$.

\sm

Let $K$ be a  metric space, $A$, $B$ be closed subsets.
Define the  {\it Hausdorff deviation} $\eta_B(A)$ as the supremum of distance between $a$ ranging in $A$ and $B$ (a number $\eta_B(A)$ is a nonnegative real or $\infty$).
The {\it Hausdorff $\infty$%
\footnote{We allow distance $=+\infty$}} on the space of closed
subset is defined by
 $$h(A,B)=\max (\eta_A(B),\eta_B(A)).$$
Its restriction to the space of compact subsets is a metric.
If  $K$ is compact then the space of its  closed subsets 
is compact.

\sm

Now we introduce the topology   on $\Mod(V)$.
We say
that $R_j$ converges to $R$ if for each
$l$ we have a convergence $B(p^l)\cap R_j \to B(p^l)\cap R$
in the sense of Hausdorff metric.
Notice that this convergence is metrized, a (non-canonical) metric is given by
$$
d(L,M)=\sum_{j=1}^\infty (2p)^{-j} h\bigl( L\cap B(p^l), M\cap B(p^l)\bigr)
.
$$

\begin{lemma}
{\rm a)} The space $\Mod(V)$ is compact with respect to the Hausdorff
topology.

\sm

{\rm b)} The space $\Lat(V)$ is a discrete dense subset
in $\Mod(V)$.
\end{lemma}

Let us prove a). Choose a convergent subsequence from
arbitrary sequence of submodules $L_j$.
First, we choose a subsequence  $L_{j_k}$ such that 
$L_{j_k}\cap B(p^0)$ converges. From the latter sequence
we choose a subsequence such that intersections
with  $B(p^1)$ converges. Etc.

\sm

{\bf\punct Analog of the radial limit.%
\label{ss:convergence-modules}}
 We need an analog of the radial limit (\ref{eq:radial}).
Say that a sequence $R_j$ of submodules
 {\it rigidly converges} to a submodule $R$  (notation $R_j\nearrow R$) if 

\sm

 {\bf(A)} for any compact subset $S\subset R$ we have $S\subset R_j$
starting some place. 

\sm

{\bf(B)} for each $\epsilon>0$, for sufficiently large $j$ the set $R_j$
is contained in the $\epsilon$-neighborhood of $R$.

\sm

{\sc Example.} Let $V=\Q_p^2$. Let $R_j=p^{-k}\Z_p e_1\oplus p^k \Z_p e_2$.
Then $R_j$ rigidly converges to a line $\Q_p e_1$.
Now let
$$
S_j= \Z_p(p^{-k} e_1+e_2)\oplus p^k \Z_p e_2
.
$$
Then $S_j$ converges to the line $\Q_p e_1$ in Hausdorff sense but not
 rigidly. 
\hfill $\boxtimes$

\sm

Evidently, we can reformulate the condition {\bf(A)} as
$$
\eta_R(R_j)\to 0
.
$$

\begin{lemma}
 The condition {\bf(B)} is equivalent to 
 $$
 \eta_{R^\lozenge}(R_j^\lozenge)\to 0
 .$$
\end{lemma}

{\sc Proof.} Let us equip $V'$ by the dual norm. Let $S$, $S_j\in V'$ and $\eta_S(S_j)\to 0$.
For small $\epsilon>0$ we have 
\begin{equation}
S_j \subset S+B(\epsilon).
\label{eq:ex-1}
\end{equation}
Passing to the duals, we get
\begin{equation}
S_j^\lozenge \supset S^\lozenge\cap B(\epsilon^{-1}).
\label{eq:ex-2}
\end{equation}
But $S^\lozenge\cap B(\epsilon^{-1})$ is an exhausting sequence of compact subsets
in $S^\lozenge$. Also, (\ref{eq:ex-2}) implies (\ref{eq:ex-1}). 
\hfill $\square$

\begin{lemma}
\label{l:1}
If $R_j\nearrow R$, then $(R_j)_\downarrow\subset R_\downarrow$  
and
$(R_j)^\uparrow  \supset R^\uparrow$ starting some $j$.
\end{lemma}

{\sc Proof.} The first claim.  For sufficiently large $k$ we have $R\subset R_\downarrow+B(p^k)$,
also $B(p^k)+B(\epsilon)=B(p^k)$ for $\epsilon\le p^{k}$.
Therefore for  a large $j$ we have
$$R_\downarrow+B(p^k)\supset R_j\supset (R_j)_\downarrow $$
But a subspace, which is contained in a tube neighborhood of  a subspace $R_\downarrow$,
is contained in $R_\downarrow$.

The second claim.
We consider a compact subset $K\subset R$ generating $R^\uparrow$ as a $\Q_p$-subspace. 
Then $(R_j)^\uparrow$ contains $K$ for sufficiently large $j$ and therefore
$(R_j)^\uparrow \supset R^\uparrow$.
\hfill $\square$

\sm

In particular, a {\it $\nearrow$-convergent sequence of linear subspaces is eventually 
constant.}

\begin{lemma}
\label{l:2}
{\rm a)}
Let $L\subset V$ be a linear subspace.
If $R_j\nearrow R$, then $(L\cap R_j)\nearrow (L\cap R)$.

\sm

{\rm b)} Let $M\subset V$ be a linear subspace, denote by
$\pi$ the natural map $V\to V/M$.
If $R_j\nearrow R$ then $\pi(R_j)\nearrow \pi(R)$. 
\end{lemma}

{\sc Proof.} a) Only condition {\bf(B)} requires a proof, i.e., 
for each $\epsilon>0$ there exists $N$ such that for $j\ge N$
$$
R_j\cap L \subset (R\cap L)+B(\epsilon)
.
$$
It is easy to shown that there is a basis $e_m$ in $\Q_p^n$ such that
$R$ has canonical form (\ref{eq:R}) and $L$ is a linear span of several basis
elements. Then for sufficiently big $N$  we have
$$
(R+ p^N\oplus \Z_p e_j)\cap L\,\,\subset \,\, (R\cap L) + p^N\oplus \Z_p e_j.
$$
Passing from the basis $e_m$ to the standard basis in $\Q_p^n$ we get 
$$
(R+B(\delta))\cap L\,\, \subset\,\, (R\cap L)+B(C\delta).
,$$
where $C=C(R,L)$ is a  constant. Now we take $\delta=\epsilon/C$
and choose number $k$, starting which $R_j \subset R+B(\delta)$.

\sm

b) follows from a) by the duality.
\hfill $\square$







\sm

{\sc Remark.}  $\nearrow$-convergence is not metrizable.
\hfill $\boxtimes$

\bigskip

{\bf B. Bruhat--Tits buildings}

\bigskip

{\bf\punct Self-dual modules.%
\label{ss:self-dual}} For details, see 
\cite{Ner-gauss}, Sections 10.6--10.7.
Consider a $p$-adic linear space $V\simeq \Q_p^{2n}$ equipped with a 
nondegenerate skew-symmetric bilinear form $B_V(\cdot,\cdot)$
(as above). We say that a subspace $L$ is {\it isotropic} 
if $B_V(v,w)=0$ for all $v$, $w\in V$.
By $\LGr(V)$ we denote the set of all maximal isotropic ({\it Lagrangian}) subspaces
in $V$ (their dimensions $=n$). 

By $L^\bot$ we denote the ortho-dual of a subspace $L$, i.e set of all vectors $w$ such that
$B_V(v,w)=0$ for all $v\in L$.

If $P$ is a submodule, denote by $P^\bbot$ the {\it dual submodule},
i.e., the set of vectors
$w$ such that $B(v,w)\in\Z_p$ for all $v\in P$.
If $P$ is a subspace, then $P^\bbot=P^\bot$.

We say that a submodule $R\subset V$ is {\it isotropic} if
$B_V(v,w)\in\Z_p$ for all $v$, $w\in R$. 

\sm

{\sc Example}. If $R$ is a linear subspace, then $R$ is isotropic in the usual sense.
On the other hand, the lattice $\Z^{2n}_p$ is isotropic (and self-dual, see below). 
\hfill $\boxtimes$
 
\sm 
 
 We say that
a submodule $R$ is {\it self-dual} if it is a  maximal isotropic
submodule in $V$. Equivalently, $P^\bbot=P$. Denote by $\LMod(V)$ the set of
all self-dual submodules in $V$,
by $\LLat(V)$ the set of all self-dual lattices. It is easy to
show that $\Sp(2n,\Q_p)$ acts on $\LLat(V)$ transitively and
$$
\LLat(V)=\Sp(2n,\Q_p)/\Sp(2n,\Z_p)
.
$$

\begin{lemma}
{\rm a)}
For any self-dual submodule $R$ 
the subspace $R_\downarrow$ is isotropic, and
$R^\uparrow$ is the  ortho-dual of $R_\downarrow$.

\sm

{\rm b)} Let $L$ ranges in the set
of   isotropic subspaces.
 Denote by $\pi:L^\bot \to L^\bot/ L$ 
the natural projection map. 
Any self-dual submodule $R$ has the form
$\pi^{-1} S$, where $S$ is 
a self-dual lattice in $L^\bot/L$.
\sm

{\rm c)} The unique $\Sp(V)$-invariant of a self-dual module $R$
is $\dim R_\downarrow$.
\end{lemma}

These statement is obvious.

\sm

Sometimes 
it is convenient to reformulate a definition of an isotropic
module. Define a bicharacter $\beta$ on
$V\times V$ by 
\begin{equation}
\beta(x,y)=\exp\bigl\{2\pi i B(x,y)\bigr\}
.
\label{eq:bich}
\end{equation}
We say that a module  $P$ is {\it isotropic} if $\beta(x,y)=1$ on $P\times P$. 


\sm

{\bf\punct Almost self-dual modules.%
\label{ss:almost}} Let $V$ and $B$ be
 same as above.
A submodule $L$ in $V$ is {\it almost self-dual} if it contains a self-dual module $M$
and $B(v,w)\in p^{-1} \Z_p$ for all $v$, $w\in L$ (see, e.g., \cite{Ner-book}, Section
10.6). Notice that 
$L/M\simeq (\Z/p \Z)^k$ with $k=0$, 1, \dots, $n$.
.
\begin{lemma}
{\rm a)}
Any almost self-dual module can be reduced by a symplectic transformation to the form
\begin{multline}
(p^{-1}\Z_p e_1\oplus \Z_p e_{n+1})
\oplus\dots\oplus
(p^{-1}\Z_p e_k\oplus \Z_p e_{n+k})
\oplus\\ \oplus
(\Z_p e_{k+1}\oplus \Z_p e_{n+k+1})
\oplus\dots\oplus
(\Z_p e_{m}\oplus \Z_p e_{n+m})
\oplus\\ \oplus
\Q_p e_{m+1}\oplus\dots\oplus \Q_p e_n
.
\end{multline}

{\rm b)} The only $\Sp(V)$-invariants of an almost self-adjoint module
$R$ are $\dim R_\downarrow$ and the number $k$ (rank of an Abelian 
group $R/S$, where $S$ is a self-dual submodule in $R$. For almost
self-dual lattices the only $\Sp(V)$-invariant is the volume of $R$, it is equal
$p^{-k}$.
\end{lemma}


{\bf\punct Graph $\Delta(V)$.} Consider a $p$-adic linear space $V$
equipped with a nondegenerate skew-symmetric bilinear form $B$  as above.
We draw an oriented graph $\Delta(V)$. Vertices are almost 
self-dual modules in $V$. 
If $R\supset R'$, then we draw an arrow from $R$ to $R'$.

\sm

If $R$, $R'$ are connected by an arrow, then $R_\downarrow=(R')_\downarrow$
and $R^\uparrow=(R')^\uparrow$.

Any pair of lattices can be connected by a (non-oriented) way.
Denote the subgraph whose vertices are all lattices by $\Delta_0(V)$.

More generally, fix an isotropic subspace
$L$ and consider the subgraph $\Delta_L(V)$ whose vertices are  almost self-dual modules
$R$ such that $R_\downarrow=L$, $R^\uparrow=L^\bot$. We get a connected
subgraph, moreover
$$
\Delta_L(V)\simeq \Delta_0(L^\bot/L)
.
$$

By definition,
$$
\Delta(V)=\bigsqcup_{\text{$L$ is isotropic subspace}} \Delta_L(V)
.
$$

If $L\subset M$, then $\Delta_M$ is contained in the closure of
$\Delta_L$ in the sense of $\nearrow$-convergence.

\sm

\begin{figure}
$$\epsfbox{build.1}$$

\caption{A reference to Subsections \ref{ss:convergence-modules}, \ref{ss:almost}.
 A subcomplex ('apartment') of the building  $\Bd(\Q_p^4)$ corresponding to 
 lattices of the form $R_1\oplus \dots\oplus R_4$, where $R_j$ is a submodule in the line $\Q_p e_j$.
 \newline
1) Vertices of the central piece of the subcomplex are almost self-dual lattices of the form
$L=p^{k_1}\Z_p e_1\oplus p^{k_2}\Z_p e_2\oplus p^{l_1}\Z_p e_3\oplus p^{l_2}\Z_p e_4$.
They are almost self-dual iff $k_1+l_1$, $k_2+l_2$ are 0 or $-1$.
\newline
2) Four boundary pieces. Each piece corresponds to almost self-dual submodules containing a line $\Q_p e_j$, where $j=1$, 2, 3, 4. 
For instance, for $j=1$  such submodules have a form
 $M=\Q_p e_1\oplus p^{m_2}\Z_p e_2\oplus  p^{l_2}\Z_p e_4$, where $m_2+l_2=0$ ,1.
A sequence of lattices $\nearrow$-converges to $M$ only if $k_1\to -\infty$ and $k_2=m_2$ starting 
some place.
\newline
3) Four extreme points correspond to Lagrangian planes spanned by pairs
of vectors
 $(e_1,e_2)$, $(e_1,e_4)$, $(e_2,e_3)$, $(e_3,e_4)$. A sequence of lattices
 $\nearrow$-converges to
$\Q_p e_1\oplus \Q_p e_4$ iff $k_1\to+\infty$, $k_2\to-\infty$.
  \label{fig:1}}
\end{figure}


{\bf\punct Bruhat--Tits buildings,%
\label{ss:buildings}} for details, see \cite{Gar}, \cite{Ner-book}.
Now we consider all $k$-plets of vertices of $\Delta(V)$ that are pairwise
connected by edges. For any such $k$-plet we draw
a $(k-1)$-simplex with given vertices and edges. Faces of a simplex
correspond to subsets of the $k$-plet. Thus we get a simplicial complex, 
denote it by
$\Bd(V)$.

Consider  the subgraph $\Delta_0$.
It can be shown that $k\le n+1$ and each simplex is contained 
in an $n$-dimensional simplex. In this way we get a structure of an
$n$-dimensional simplicial complex, it is called a {\it Bruhat--Tits building.}
We denote it by $\Bd_0(V)$.

For a subgraph $\Delta_L$ we get a simplicial complex 
complex $\Bd_L(V)$ isomorphic $\Bd(L^\bot/L)$. 

Below we use term '{\it distinguished boundary of a building}' for the
Lagrangian Grassmannian, this is an counterpart of Shilov boundary.


\bigskip

{\bf C. Relations and Nazarov category}

\bigskip


{\bf\punct Relations.%
\label{ss:relations}} Let $V$, $W$ be linear spaces.
We say that a {\it relation} $P:V\tto W$ is a submodule in
$V\oplus W$.

\sm

{\sc Example.} Let $A:V\to W$ be a linear operator. Then its graph is a relation.
\hfill$\boxtimes$

\sm

Let $P:V\tto W$, $Q:W\tto Y$ be relations.
We define their product $S=QP: V\tto Y$ as the set of all 
$v\oplus y\in V\oplus Y$ for which
 there exists $w\in W$ such that
$v\oplus w\in P$, $w\oplus y\in Q$.

\sm

For a relation $P:V\tto W$ we define its {\it kernel} $\ker P\subset V$ as
$$
\ker P= P\cap (V\oplus 0)
,$$
 the  {\it indefiniteness} $\indef P\subset W$,
$$
\indef P= P\cap (0\oplus W)
,$$
the {\it domain of definiteness}
$$
\dom P= \text{projection of $P$ to $V$}
,$$
and 
the {\it image}
$$
\im P= \text{projection of $P$ to $W$}
.
$$
We  define the {\it pseudo-inverse} relation
$P^\square:W\tto V$ being the same submodule in $W\oplus V\simeq V\oplus W$.
Evidently,
$$
(PQ)^\square=Q^\square P^\square
.
$$


{\bf\punct The definition of product. A reformulation.%
\label{eq:def-product-2}}
We keep the same notation.
Consider the space $\cZ:=V\oplus W\oplus W \oplus Y$
and following submodules of $\cZ$:

\sm

--- the subspace $\cH$ consisting of vectors $v\oplus w\oplus w\oplus y$;

\sm

--- the subspace $\cA$ consisting of vectors $0\oplus w\oplus w\oplus 0$;

\sm

--- the submodule $P\oplus Q\subset (V\oplus W)\oplus (W\oplus Y)$.

\sm

Then we do the following operations:

\sm

--- take the intersection $R=\cH\cap (P\oplus Q)$;

\sm

--- take the map $\theta:\cH\to \cH/\cA\simeq V\oplus Y$.

\sm 

Then $QP=\theta(R)$.

\sm


{\bf\punct  Action on $\Mod (V)$.%
\label{ss:action-submodules}}
Let $P:V\tto W$ be a  relation, $T$ be a submodule in $V$.
We define the submodule $PT\subset W$ as the set of 
$w\in W$ such that there is $v\in T$ satisfying $v\oplus w\in P$.

\sm

{\sc Remark.}
We can consider a submodule $T\subset V$ as a relation  $0\tto V$.
Therefore we can regard $PT:0\tto W$ as the product of relations $T: 0\tto V$
and $Q: V\tto W$. \hfill $\boxtimes$

\sm


{\bf \punct The Nazarov category.%
\label{ss:nazarov}}
For a pair $V$, $W$ of symplectic linear spaces
we define a skew-symmetric bilinear form $B^\ominus$ on $V\oplus W$ by
$$
B^\ominus(v\oplus w,v'\oplus w')
=B_V(v,v')-B_W(w,w') 
.
$$
Denote by

\sm

--- $\ov\Naz(V,W)$ the set of all self-dual submodules of $V\oplus W$;

\sm

--- $\Naz(V,W)$ the set of $P\in\ov\Naz(V,W)$ such that $\ker P$ and $\indef P$ are compact.

\begin{theorem}
\label{th:nazarov-0}
Let $P\in\ov\Naz(V,W)$, let $T$ be a self-dual submodule
in 
 $V$. Then the submodule  $PT\subset W$ is self-dual.
\end{theorem}

In \cite{Ner-gauss}, Theorem 10.7.2,  the same statement
is established under slightly stronger condition
 $P\in\Naz(V,W)$. In fact, a proof remains valid
 for  
$P\in\ov\Naz(V,W)$.

\begin{theorem}
\label{th:nazarov-1}
{\rm a)} If $P\in\Naz(V,W)$, $Q\in \Naz(W,Y)$, then $QP\in \Naz(V,Y)$.

\sm

{\rm b)} If $P\in\ov\Naz(V,W)$, $Q\in \ov\Naz(W,Y)$, then $QP\in \ov\Naz(V,Y)$.

\sm

{\rm c)} If $P\in\Naz(V,W)$, $Q\in \Naz(W,Y)$ are lattices, then $QP$ is a lattice.
\end{theorem} 

The statement a) was proved in Nazarov \cite{Naz} (see also \cite{Ner-gauss}, Section 10.7), c) is obvious.
The statement b) is a corollary of Theorem
 \ref{th:nazarov-0}, see
\cite{Ner-gauss}, Subsection 10.7.4.

\sm

Thus we get two similar categories%
\footnote{The Nazarov category is an analog of Krein--Shmulian type categories,
see \cite{Ner-book}, \cite{Ner-gauss}},
 $\Naz$ and $\ov\Naz$. {\it The group of automorphisms of an object $V$ is
the symplectic group $\Sp(V,\Q_p)$} (for both categories),
an operator $V\to V$ is symplectic iff its graph is isotropic with respect to the form
$B^\ominus$.

\sm
 
For $P\in\ov\Naz(V,W)$, we have
\begin{align*}
&(\ker P)^\bbot=\dom P,\qquad & &(\indef P)^\bbot=\im P
\\
&\bigl((\ker P)_\downarrow\bigr)^\bot=(\dom P)^\uparrow,
\qquad
&&\bigl((\indef P)_\downarrow\bigr)^\bbot=(\im P)^\uparrow,
\end{align*}
 
{\bf\punct Action of the Nazarov category on buildings.%
\label{ss:naz-action}}

\begin{proposition}
\label{prop:almost-almost}
{\rm a)} Let $P\in \Naz(V,W)$, let $T$ be an almost-self-dual lattice.
Then $PT\subset W$ is an almost self-dual lattice.

{\rm b)} Let $P\in \ov\Naz(V,W)$, let $T$ be an almost-self-dual submodule.
Then $PT\subset W$ is an almost self-dual submodule.
\end{proposition}

 The statement a) is \cite{Ner-gauss}, Proposition 10.7.5, a proof
 remains to be valid for the statement b) also.
 
 Now, let $\Xi$, $\Sigma$ be simplicial complexes,
 let $\Ver(\Xi)$, $\Ver(\Sigma)$ be their sets of vertices.
 We say, that a map%
 \footnote{generally, non-injective.}
  $\Ver(\Xi)\to\Ver(\Sigma)$
 is {\it simplicial}, if for any simplex $\Delta\subset \Xi$
 images of its vertices are are contained in one simplex of $\Sigma$.
 Notice, that we can extend a simplicial map to a map of
 complexes $\Xi\to\Sigma$ assuming that a map is affine 
 on each face.

The following statement is a corollary of Proposition \ref{prop:almost-almost}.

\begin{theorem}
\label{th:naz-action}
{\rm a)} A morphism $P\in \Naz(V,W)$ induces  simplicial map
$$\Bd(V)\to \Bd(W)$$
sending
 $$\Bd_0(V)\to \Bd_0(W) .$$
{\rm b)} A morphism $P\in \ov\Naz(V,W)$ induces a simplicial map
$\Bd(V)\to \Bd(V)$, sending $\Bd(V)$ to
\begin{equation}
\Bd(V)
\to
\bigsqcup_{\begin{matrix}\text{\rm\small $M$ is isotropic subspace in $W$}
\\ M\supset
\indef(P)_\downarrow
\end{matrix}
} \Bd[M^\bot/M]
.
\label{eq:bigsqcup}
\end{equation}
\end{theorem}

{\sc Remark.}
The map $T\to PT$ is contractive in an essentially stronger
sense, see \cite{Ner-comp}.
\hfill$\boxtimes$

\begin{theorem}
\label{th:continuity0}
Let $P\in\ov\Naz(V,W)$. The the induced map 
$\Bd(V)\to\Bd(W)$
is $\nearrow$-continuous, i.e., for a convergent
sequence
$T_j \nearrow T$ of almost self-dual modules,
we have $PT_j \nearrow PT$.
\end{theorem}  

{\sc Proof.} We evaluate $PT_j$ according procedure
described in Subsection \ref{eq:def-product-2}.
By Lemma \ref{l:2}, both steps of the evaluation are continuous.
\hfill $\square$

\bigskip

{\bf D. Weil representation}

\bigskip

The Weil representation is used below only in Section 7.

\sm


{\bf\punct Extended Nazarov category.%
\label{ss:extended-nazarov}} Now we add to the Nazarov category an
infinite-dimensional object $V_{2\infty}$.
This is the space of vectors
$$
(x_1^+,x_2^+,\dots, x_1^-,x_2^-,\dots),\qquad \text{where $x_j^\pm\in\Q_p$ and $x_j^\pm\in\Z_p$
for almost all $j$}
.
$$
Notice that $V_{2\infty}$ is not a $\Q_p$-linear space but is
a $\Z_p$-module.

We introduce a bicharacter $\beta(\cdot,\cdot)$ on $V_{2\infty}\oplus V_{2\infty}$ by
$$
\beta(x,y)= \exp\Big[2\pi i\sum_{j=1}^\infty (x_j^+ y_j^- - x_j^- y_j^+)\Bigr]:=
\prod_{j=1}^\infty
\exp\bigl\{2\pi i(x_j^+ y_j^-  -x_j^- y_j^+)\bigr\}
.
$$
Notice that almost all factors of the product equal to 1. The sum in square brackets
defining a symplectic form
is not well defined, more precisely it is well defined modulo $\Z_p$.

Objects of the {\it extended Nazarov category} $\bNaz$ are

\sm

---  finite-dimensional spaces
$V$ equipped with skew-symmetric non-degenerate bilinear forms $B_V$
and with the corresponding bicharacters $\beta_V$, see (\ref{eq:bich});

\sm

--- the space $V_{2\infty}$.

\sm

Let $V$, $W$ be two objects. We equip their direct sum
with a bicharacter
$$
\beta_{V\oplus W}(v\oplus w,v'\oplus w')=\frac{\beta_V(v,v')}{\beta_W(w,w')}
.
$$
A {\it morphism of the category $\bNaz$} is a  self-dual submodule
 $P\subset V\oplus W$ such that $\ker P$ and $\indef P$ are compact.
 
 Group $\mathbf{Sp}(2\infty,\Q_p)$ 
 of automorphisms of $V_{2\infty}$ consists of $2\infty\times 2\infty$ matrices
 $r=\begin{pmatrix}a&b\\c&d\end{pmatrix}$ such that
 
 \sm
 
  --- all but a finite number of matrix elements are integer;
  
\sm  
  
  --- matrix elements $a_{ij}$, $b_{ij}$, $c_{ij}$, $d_{ij}$ tend to 0 as $i\to\infty$ for fixed $j$; also 
  they tend to 0 as $j\to\infty$ for fixed $i$;

\sm

 --- matrices $r$ are symplectic in the usual sense, 
 $$
 r^t\begin{pmatrix}0&1\\-1&0\end{pmatrix} r=\begin{pmatrix}0&1\\-1&0\end{pmatrix}=
r \begin{pmatrix}0&1\\-1&0\end{pmatrix} r^t
.
$$


{\bf\punct Heisenberg groups.%
\label{ss:heis}} For the sake of simplicity, set $p>2$.
Denote by $\T_p\subset \C^\times$ the group of complex roots
of unity of  degrees $p$,  $p^2$, $p^3$,\dots.
Let $V$ be an object of the extended Nazarov category. We define the Heisenberg group
$\Heis(V)$ as a central extension of the Abelian group $V$ by $\T_p$ in the following way.
As a set, $\Heis(V)\simeq V\times \T_p$. The multiplication is given by
$$
(v,\lambda)\cdot (w,\mu)=
\bigl(v+w, \lambda\mu \cdot \beta_V(v,w)\bigr)
.
$$

Decompose $V=V_+\oplus V_-$ as in (\ref{eq:V-pm}).
For a finite dimensional $V$ we define a unitary representation
$\Psi$
of $\Heis(V)$ in $L^2(\Q_p^n)$ by the formula
\begin{equation}
\Psi(v^+\oplus v^-,\lambda) f(x)=
\lambda f(x+v^+)\,\exp\Bigl\{2\pi i\Bigl(\sum v_j^+ x_j+\frac 12 \sum v_j^+ v_j^-\Bigr)\Bigr\}
\label{eq:heis}
.
\end{equation}

Next, consider the space $\cE_\infty$ consisting of sequences
 $z=(z_1, z_2,\dots)$ such that
$|z_j|\le 1$ for all but a finite number of $j$. This space
is an Abelian locally compact group, it admits a Haar measure.
On the open subgroup $\Z_p^\infty\subset \cE_\infty$, the
Haar measure  is a product of probability 
Haar measures on $\Z_p$. The whole space $\cE_\infty$ is a countable disjoint union
of sets $u+\Z_p^\infty$. 

We define the representation of the group $\Heis(V_{2\infty})$ in 
$L^2(\cE_\infty)$ by the same formula (\ref{eq:heis}).

\sm


{\bf\punct The  Weil representation of the Nazarov category. Formal definition.%
\label{ss:weil}}
See \cite{Naz}, \cite{NNO}, for finite-dimensional case, see \cite{Ner-gauss}, Chapter 11.

\begin{theorem}
\label{th:nazarov-2}
 For a $2n$-dimensional object of the category
$\bNaz$ we assign the Hilbert space $\cH(V):=L^2(\Q_p^n)$. For the object $V_{2\infty}$, we assign
the Hilbert space $\cH(V_{2\infty}):=L^2(\cE_\infty)$. 

\sm

{\rm a)} Let $V$, $W$ be objects of $\bNaz$. Let $P:V\tto W$ be 
a morphism of category $\bNaz$. Then there is a unique up to a scalar factor bounded operator
$$\We(P):\cH(V)\to\cH(W)$$
 such that
$$
\Psi(w,1)\We(P)=\We(P)\Psi(v,1)
\quad\text{for all $v\oplus w\in P$} 
.
$$

\sm

{\rm b)} Let $V$, $W$, $Y$ be objects of $\bNaz$. Let $P:V\tto W$, $Q: W\tto Y$ be morphisms
of $\bNaz$. Then
$$
\We(Q)\We(P)=s\cdot \We(QP)
,$$
where $s=s(P,Q)\in\C^\times$ is a nonzero scalar. In other words, we get a  {\rm projective
representation} of the category $\bNaz$. Also,
$$\We(P^\square)=t\cdot\We(P)^*,\qquad t\in\C^\times.$$
\end{theorem}

For symplectic groups $\Sp(2n,\Q_p)=\mathrm{Aut}(\Q_p^{2n})$ the representation $\We(g)$
 coincides with the Weil representation. 

\sm

{\bf \punct
Explicit formulas for operators  for some  morphisms.}
 
 1) Let
$V=W$ and $P$ be a graph of a symplectic  operator. There are simple formulas
for some special symplectic matrices:
\begin{align}
\We\begin{pmatrix}
A&0\\0&A^{t-1}
\end{pmatrix}
f(z)&=|\det A|^{1/2} f(zA)
;
\label{eq:weil-gl}
\\
\We\begin{pmatrix}
1&B\\0&1
\end{pmatrix}
f(z)&= \exp\{\pi i zBz^t\}
;
\nonumber
\\
\We\begin{pmatrix}
0&1\\-1&0
\end{pmatrix}
f(z)&=\int_{\Q_p^n} f(x)\exp\{2\pi i xz^t\}\,dx
.
\nonumber
\end{align} 
Any element of $\Sp(2n,\Q_p)$ can be represented as a product of matrices
of such forms, this allows to write an explicit formula for $\We(g)$ for any element
 $g\in\Sp(2n,\Q_p)$. 

\sm

Denote by $I(x)$ the function on $\Q_p$ defined by 
$$
I(x)=\begin{cases}
1,\qquad |x|\le 1;
\\
0,\qquad\text{otherwise}.
\end{cases}
$$

Next, we need some special non-invertible morphisms.

\sm

2) Let $V=\Q_p^{2n}$, $W=V\oplus Y$, where $Y=\Q_p^{2n}$ or $V_{2\infty}$.
Denote by $Y(\Z_p)$ the lattice $\Z_p^{2n}$ or $\Z_p^{2\infty}$
respectively.
 Denote by
$$
\lambda^V_W:V\tto W
$$
 the direct sum of the graph $\graph(1_V)$
  of the unit operator $1_V:V\to V$ and the lattice
$Y(\Z_p)\subset Y$. Then
$$
\We(\lambda_W^V)\, f(v_1,\dots,v_n,y_1, y_2,\dots)=f(v_1,\dots,v_n)\,I(y_1) I(y_2) \dots
$$

3) Preserving the previous notation denote by 
$$
\theta^V_W:W\tto W
$$
 the direct sum
$$
\graph(1_V)\oplus(Y(\Z_p)\oplus Y(\Z_p))\,\subset\, (V\oplus V)\oplus ( Y\oplus Y)
.
$$
Then
\begin{equation}
\theta_W^V=\lambda_W^V\left(\lambda_W^V\right)^* ,\qquad 
\left(\theta_W^V\right)^2=\theta_W^V,\qquad \left(\lambda_W^V\right)^*\lambda_W^V=1_V.
\label{eq:lambda-theta}
\end{equation}
The operator
$\We(\theta_W^V)$ is the orthogonal  projection to the space of functions of the form
$$
f(v_1,\dots,v_n)\,I(y_1)I(y_2) \dots 
$$

{\bf\punct General case.}
Any  morphism of the category $\bNaz$ can be represented as a product of morphisms of the types
described above.
Moreover, for finite dimensional $V$, $W$, any $P:V\tto W$ can be represented as
$$
P=
(\lambda_Z^W)^* \cdot g\cdot \lambda_Z^V, \qquad  
g\in \Sp(Z)
,$$
where $Z$ is sufficiently large ($\dim Z\ge 2\max(\dim V,\dim W)$).
In fact, the same decomposition holds for morphisms $Q:V_{2\infty} \to V_{2\infty}$,
any $Q$ can be represented as
$$
Q=
\theta^{V_{2\infty}}_{V_{2\infty}\oplus V_{2\infty}}
\cdot
g
\cdot
\theta^{V_{2\infty}}_{V_{2\infty}\oplus V_{2\infty}},
\qquad g\in \Sp(V_{2\infty}\oplus V_{2\infty}).
$$


\section{Characteristic function}

\COUNTERS

Here we define characteristic functions of double cosets
$\bfK\setminus \bfG/\bfK$ and formulate several theorems. Proofs are in the next
section.

\sm


{\bf\punct Construction.%
}
Consider the group
$$
\GL(\alpha+k\infty,\Q_p)
:=\lim_{j\to\infty} \GL(\alpha+kj,\Q_p).
$$
Let $g\in \GL(\alpha+k\infty,\Q_p)$ actually be contained in $\GL(\alpha+km,\Q_p)$,
\begin{equation}
g
=
\begin{pmatrix}
a&b_1&\dots& b_k\\
c_1&d_{11}&\dots&d _{1k}\\
\vdots&\vdots & \ddots&\vdots\\
c_k&d_{k1}&\dots& d_{kk}
\end{pmatrix}
\in\GL(\alpha+km,\Q_p)
\label{eq:g}
.
\end{equation}
We write the following equation
(this is an analog of (\ref{eq:char-uo-infty}), the analogy is important)
\begin{equation}
\begin{pmatrix}
v^+\\
y_1^+\\
\vdots\\
y_k^+\\
v^-\\
y_1^-\\
\vdots\\
y_k^-\\
\end{pmatrix}
=
\begin{pmatrix}
\begin{matrix}
a&b_1&\dots& b_k\\
c_1&d_{11}&\dots&d _{1k}\\
\vdots&\vdots & \ddots&\vdots\\
c_k&d_{k1}&\dots& d_{kk}
\end{matrix}
&
\begin{matrix}
0_{\phantom{1}}&0_{\phantom{11}}&\dots&0_{\phantom{1k}}\\
0_{\phantom{1}}&0_{\phantom{11}}&\dots&0_{\phantom{1k}}\\
\vdots_{\phantom{1}}&\vdots_{\phantom{1k}} & \ddots_{\phantom{1}}&\vdots_{\phantom{1k}}\\
0_{\phantom{1}}&0_{\phantom{1k}}&\dots&0_{\phantom{1k}}
\end{matrix}
\\
\begin{matrix}
0_{\phantom{1}}&0_{\phantom{11}}&\dots&0_{\phantom{1k}}\\
0_{\phantom{1}}&0_{\phantom{11}}&\dots&0_{\phantom{1k}}\\
\vdots_{\phantom{1}}&\vdots_{\phantom{1k}} & \ddots_{\phantom{1}}&\vdots_{\phantom{1k}}\\
0_{\phantom{1}}&0_{\phantom{1k}}&\dots&0_{\phantom{1k}}
\end{matrix}
&
\,\,\,\,
\begin{pmatrix}
a&b_1&\dots& b_k\\
c_1&d_{11}&\dots&d _{1k}\\
\vdots&\vdots & \ddots&\vdots\\
c_k&d_{k1}&\dots& d_{kk}
\end{pmatrix}^{t-1}
\end{pmatrix}
\begin{pmatrix}
u^+\\
x_1^+\\
\vdots\\
x_k^+\\
u^-\\
x_1^-\\
\vdots\\
x_k^-\\
\end{pmatrix}
\label{eq:main-1}
.
\end{equation}
Here
$u^\pm$,  $v^\pm\in \Q_p^\alpha$
 and $x^\pm_j$, $y^\pm_j\in
 \Q_p^m$. 
 
Before the exploring of this 
 identity as (\ref{eq:char-uo-infty}), 
 we need some preparations.
 
 \sm
 
Define 3 spaces, $\cV$, $\cH$, $\ell_m$: 

\sm

1) Denote $\cV:=\Q_p^\alpha\oplus \Q_p^\alpha$. We regard
$u=u^+\oplus u^-$, $v=v^+\oplus v^-$ as elements of $\cV$.
Equip $\cV$ with the standard skew-symmetric bilinear  form 
$\begin{pmatrix} 0&1\\-1&0\end{pmatrix}$.

\sm

2) Denote 
\begin{equation}
\cH:=\cH^+\oplus \cH^-=\Q_p^k\oplus \Q_p^k
\label{eq:H}
\end{equation}
and equip this space with the standard skew-symmetric bilinear form.

\sm

3) Denote by $\ell_m$ the space $\Q_p^m$ equipped with the standard symmetric bilinear form
$$(z,w)=\sum z_j w_j.$$
 We regard $x^\pm_j$, $y^\pm_j$ as elements of this space.

\sm

Consider  the tensor product $\cH\otimes_{\Q_p} \ell_m$,
vectors 
$$
\begin{pmatrix}x_1^+& \dots& x_k^+& x_1^-& \dots& x_k^-\end{pmatrix},
\quad
\begin{pmatrix}y_1^+& \dots& y_k^+& y_1^-& \dots& y_k^-\end{pmatrix}
$$
are regarded as elements of $\cH\otimes \ell_m$. 
We equip $\cH\otimes \ell_m$ with the tensor product of bilinear 
forms, this form is a skew-symmetric with matrix%
\footnote{A tensor product of a symmetric and
a skew-symmetric bilinear forms is a skew-symmetric
bilinear form.} 
{\small
$$
\begin{pmatrix}
0&\dots&0& 1_m&\dots&0\\
\vdots&\ddots&\vdots &\vdots&\ddots&\vdots &\\
0&\dots&0& 0&\dots&1_m\\
-1_m&\dots&0& 0&\dots&0\\
\vdots&\ddots&\vdots &\vdots&\ddots&\vdots &\\
0&\dots&-1_m& 0&\dots&0\\
\end{pmatrix}
.$$}

Thus the operator in (\ref{eq:main-1}) is an operator
$$
\cV\oplus (\cH\otimes \ell_m) 
\quad \to\quad
\cV\oplus (\cH\otimes \ell_m) 
$$
We equip the spaces $\cV\oplus (\cH\otimes \ell_m) $
with a skew-symmetric bilinear form that is a direct sum
of forms in
 $\cV$ and $\cH\otimes \ell_m$.
 The matrix of this form is
 $$
\begin{pmatrix}
0&1_\alpha&0&0\\
-1_{\alpha}&0&0&0\\
0&0&0&1_{km}
\\0&0&-1_{km}&0
\end{pmatrix}
$$   
Evidently,
 {\it operators {\rm(\ref{eq:main-1})} preserve this form},
 i.e., they are contained in $\Sp\bigl(2(\alpha+km),\Q_p\bigr)$.

\sm

Now we start a description of characteristic functions.

\sm

For any self-dual module $Q\subset \cH$
we consider the self-dual module 
$$
Q\otimes_{\Z_p} \Z_p^m \subset  \cH\otimes \ell_m
.
$$
Notice, that $Q\otimes\Z_p^m$ is a direct sum of $m$ copies of $Q$.

\sm

\begin{definition} Fix $g$.
Fix self-dual submodules $Q$, $T \subset \cH$.
We define a relation 
$$\chi_g(Q,T):\cV\tto \cV$$
 as the
set of all $u\oplus v\in \cV\oplus \cV$ for which 
there exist $x\in Q\otimes \Z_p^m$, $y\in T\otimes \Z_p^m$
such that {\rm (\ref{eq:main-1})} holds.
\end{definition}

{\bf\punct An auxiliary definition.} 

\begin{definition}
We say that some property of a double coset {\rm holds in a general position}
if for any sufficiently large $m$ the set of points
$g\in\GL(\alpha+km,\Q_p)$, where the property
does not hold, is a proper algebraic subvariety in $\GL(\alpha+km,\Q_p)$.
\end{definition}

\sm


{\bf \punct Basic properties of characteristic functions.}

\begin{lemma}
\label{l:choice}
$\chi_g(Q,T)$ does not depend on a choice of $m$.
\end{lemma}

\begin{theorem}
\label{th:independence}
If $g_1$, $g_2$ are contained in the same double 
coset $\bfK\setminus \bfG/\bfK$, then  $\chi_{g_1}(Q,T)=\chi_{g_2}(Q,T)$.
\end{theorem}

Thus, for any double coset $\frg\in \bfK\setminus \bfG/\bfK$ we get a well-defined
map 
$$
\chi_\frg:\,\,
\LMod(\cH)\times \LMod(\cH) \quad\to \quad \Bigl\{\text{space of relations $\cV\tto \cV$}\Bigr\}
.
$$
Therefore, we can write
$$
\chi_{\frg}(Q,T),\qquad\text{where $\frg\in \bfK\setminus \bfG/\bfK$}
.
$$
We say that $\chi_\frg(\cdot,\cdot)$ is the {\it characteristic function}
of the double coset $\frg$. 

\begin{theorem}
\label{th:self-duality}
$\chi_\frg(Q,T)\in \ov\Naz(\cV,\cV)$.
\end{theorem}

\begin{theorem}
\label{th:multiplication}
 The following identity holds
$$\chi_{\frg\star \frh}(Q,T)= \chi_\frg(Q,T)\, \chi_\frh(Q,T),$$
in the right-hand side we have a product of relations
\end{theorem}


{\bf\punct Refinement of Theorem \ref{th:self-duality}.%
\label{ss:refinement}}
Fix a double coset  $\frg$. Substituting $x^\pm=0$, $y^\pm=0$ 
to the equation (\ref{eq:main-1}), we get an equation
for $u\oplus v\in \cV\oplus \cV$.
The explicit form  (see equation (\ref{eq:main-2})) is 
\begin{equation}
\begin{cases}
v^+= au^+\\
0\,\,\,\,=c_j u^+,\qquad \text{for all $j$}
\\
u^-= a^tv^-\\
0\,\,\,\,=b^t_j v^-,\qquad \text{for all $j$}
\end{cases}
\label{eq:Lambda}
\end{equation}

 Denote by $\Lambda(\frg)\subset \cV\oplus \cV$
the linear subspace of solutions of this system.

Notice that
$$
\ker \Lambda(\frg)=0,\qquad \indef\Lambda(\frg)=0
$$
(since $g$ is an invertible matrix).

For $\frg$ being in a general position $\Lambda(\frg)=0$. 

\begin{proposition}
\label{pr:Lambda}
{\rm a)} For any self-dual $Q$, $T\in \LMod(\cH)$,
$$
\chi_\frg(Q,T)_\downarrow\supset \Lambda(\frg),
\qquad
\chi_\frg(Q,T)^\uparrow\subset \Lambda(\frg)^\bot.
$$

{\rm b)} If $Q$, $T$ are self-dual lattices, then
$$
\chi_\frg(Q,T)_\downarrow= \Lambda(\frg),
\qquad
\chi_\frg(Q,T)^\uparrow= \Lambda(\frg)^\bot.
$$
\end{proposition}

\begin{corollary}
 For $\frg$ being in a general position,
we get a map
$$
\LLat(\cH)\times\LLat(\cH)\,\to \, \LLat(\cV\oplus \cV)
.
$$
\end{corollary}


{\bf\punct Values of characteristic functions on the distinguished boundary.}

\begin{theorem}
\label{th:absolute}
Let  $Q$, $T$ range in the   Lagrangian Grassmannian $\LGr(\cH)$. Then

\sm

{\rm a)}  
$\chi_\frg(Q,T)$ is a Lagrangian subspace in $\cV\oplus \cV$.

\sm

{\rm b)} The map 
$$
\chi_\frg:\,
\LGr(\cH)\times \LGr(\cH)\to \LGr(\cV\oplus \cV)
$$
is rational.

\sm

{\rm c)} For $\frg$ being in a general
position,  $\chi_\frg(Q,T)\in \Sp(\cV,\Q_p)$ a.s. on $\LGr(\cH)\times \LGr(\cH)$.
\end{theorem}

A precise description of the subset of $\bfK\setminus \bfG/\bfK$, where the last property holds, is given below
in Subsection \ref{ss:determinant}. 

There is a more exotic statement in the same spirit. 

\begin{proposition}
\label{pr:exotic}
For all $\frg$ for almost all $(Q,T)\in \LGr(\cH)\times\LGr(\cH)$,
the condition $(u^+\oplus u^-)\oplus (v^+\oplus v^-)\in \chi_\frg(Q,T)$
can be written as an  equation
$$
\begin{pmatrix}
v^+\\
u^-
\end{pmatrix}
=
Z(Q,T)
\begin{pmatrix}
v^-\\
u^+
\end{pmatrix}
$$
there  $Z(Q,T)$ is a symmetric matrix.
\end{proposition}

Point out that this can done {\it for all} 
 $\frg$.

\begin{proposition}
\label{pr:another-quotient}
Let
$$\frg_1,\, \frg_2\in \bfK\setminus\bfG/\bfK=
\OO(\infty,\Z_p)\setminus\bfG/\OO(\infty,\Z_p)$$
 be contained in the same double coset
$$
\OO(\infty,\Q_p)\setminus\bfG/\OO(\infty,\Q_p)
,
$$ 
then the restrictions of $\chi_{\frg_1}$ and $\chi_{\frg_2}$
to $\LGr(\cH)\times \LGr(\cH)$ coincide.
\end{proposition}


\sm

{\bf\punct Extension of characteristic function  to buildings.%
\label{ss:maps-buildings}}
Next,  consider two almost self-dual submodules
$Q$, $T$  and apply to them the definition of characteristic function
$Q$, $T$. 

\begin{figure}
$$\epsfbox{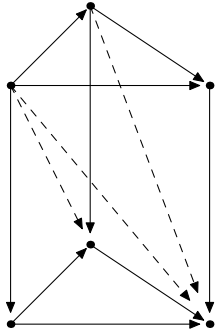}$$
\caption{A reference to Subsection  \ref{ss:maps-buildings}. A product of two simplices and additional arrows.%
\label{fig:2}}
\end{figure}

\begin{proposition}
\label{pr:almost-self-dual}
If $Q$, $T$ are almost self-dual modules, then $\chi_\frg(Q,T)$ is almost 
self-dual.
\end{proposition}

Now we construct an oriented graph $\Delta(\cH\Join \cH)$. Vertices
are ordered pairs $(Q,T)$ of almost self-dual submodules in $\cH$. We draw
an arrow from $(Q,T)$ to $(Q',T')$ if $Q\supset Q'$, $T\supset T'$.

Consider the product
of simplicial complexes $\Bd(\cH) \times \Bd(\cH)$.
It is polyhedral complex, whose cells are products of simplices.
Two vertices (of this complex) $(Q,T)$ and $(Q',T')$ are connected by an arrow
if $Q\supset Q'$ and $T=T'$ or $Q= Q'$ and $T\supset T'$.
However, our rule from the previous paragraph
 produces more arrows, this provides a simplicial 
partition of each product of simplices (see, e.g., \cite{Hat}, Section 3.B).
Finally, we get 
a $2k$-dimensional simplicial complex $\Bd(\cH\Join \cH )$
(it  also is a subcomplex of the 
 complex
$\Bd(\cH\oplus \cH)$).

Let $\Phi$, $\Psi$ be two oriented graphs, assume that number of edges connecting
any pair of vertices is $\le 1$. We say that
a map $\sigma:\Ver(\Phi)\to \Ver(\Psi)$ is a {\it morphism of graphs}
if for any arrow $a\to b$ in $\Phi$ we have 
$\sigma(a)= \sigma(b)$ or there is an arrow $\sigma(a)\to\sigma(b)$.

\begin{figure}
$$\epsfbox{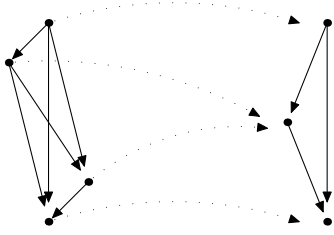}$$
\caption{A reference to Subsection \ref{ss:maps-buildings}. A morphism of oriented
graphs\label{fig:3}}
\end{figure}

\begin{theorem}
\label{th:morphisms-graphs}
A characteristic function $\chi_\frg$
is a morphism of oriented graphs 
\begin{equation}
\Delta(\cH\Join \cH)\to\Delta(\cV\oplus \cV).
\label{eq:dd}
\end{equation}
\end{theorem}



{\bf\punct Continuity.}

\begin{theorem}
\label{th:continuity}
Let $Q_j$, $Q$, $T_j$, $T$ be almost self-dual modules.
If $Q_j\nearrow Q$, $T_j\nearrow T$, then
$$
\chi_\frg(Q_j,T_j)\nearrow \chi_\frg(Q,T)
.
$$
\end{theorem}

Notice that characteristic function can be discontinuous
with respect to the Hausdorff convergence. Moreover, the restriction 
of $\chi_\frg$ to $\LGr(\cH)\times\LGr(\cH)$ can be discontinuous in the topology
of Grassmannian.


\sm

{\bf \punct Involution.}

\begin{proposition}
\label{pr:involution}
If $u\oplus v\in \chi_\frg(Q,T)$, then $v\oplus u\in \chi_{\frg^*}(T,Q)$.
\end{proposition}

\sm


{\bf\punct Additional symmetry.} For a nonzero $\lambda\in \Q_p^\times=\Q_p$,
we define an operator $M(\lambda)$ in $\cH$ given by $\begin{pmatrix}\lambda&0\\0&\lambda^{-1}
\end{pmatrix}$,
by the same symbol we denote the operator $\begin{pmatrix}\lambda&0\\0&\lambda^{-1}\end{pmatrix}$
in the space $\cV$.

\begin{theorem}
\label{th:symmetry}
$$\chi_\frg\bigl(M(\lambda)Q, M(\lambda) T\bigr)=M(\lambda^{-1}) \chi_\frg(Q,T)M(\lambda).$$
\end{theorem}


\sm


{\bf \punct Remark. Another semigroup of double cosets.}
Consider the group $\wt \bfG=\Sp(2\alpha+2k\infty, \Q_p)$ of symplectic
matrices $\begin{pmatrix}a&b\\c&d\end{pmatrix}$
of size $(\alpha+k\infty)+(\alpha+k\infty)$, $\wt \bfG\supset \bfG$. Consider its subgroup   
$\bfG= \GL(\alpha+k\infty,\Q_p)$
consisting of matrices $\begin{pmatrix}g&0\\0&g^{t-1}\end{pmatrix}$,
consider the same 
$\bfK=\OO(\infty,\Z_p)\subset \GL(\alpha+k\infty,\Q_p)$.
Consider the semigroup of double cosets
$\bfK\setminus \wt \bfG/\bfK$, the multiplication is determined as in Theorem
\ref{th:umnozhenie}.

 We define characteristic function $\chi_{\wt \frg}(Q,T)$ in the same way,
in formula (\ref{eq:main-1}) instead the  matrix
 $\begin{pmatrix}g&0\\0&g^{t-1}\end{pmatrix}$ we write a symplectic 
matrix $\begin{pmatrix}a&b\\c&d\end{pmatrix}\in \Sp(2\alpha+2k\infty, \Q_p)$. 

\begin{theorem}
All the statements of this section hold for  $\chi_{\wt \frg}(Q,T)$ 
except Theorem {\rm \ref{th:symmetry}} and Proposition \ref{pr:exotic}%
\footnote{the system (\ref{eq:Lambda}) also must be modified.}. 
\end{theorem}



\section{Proofs}

\COUNTERS

{\bf\punct Independence of representatives.} To shorten expressions, set
$k=2$. Let $h\in\OO(m,\Z_p)$, let $\frI(h)$ be given by (\ref{eq:tau}). Then characteristic
function of $g\frI(h)$ is determined by
$$
\begin{pmatrix}
v^+\\
y_1^+\\
y_2^+\\
v^-\\
y_1^-\\
y_2^-\\
\end{pmatrix}
=
\begin{pmatrix}
\begin{matrix}
a&b_1h& b_2h\\
c_1&d_{11}h&d _{12}h\\
c_2&d_{21}h& d_{22}h
\end{matrix}
&
\begin{matrix}
0_{\phantom{1}}&0_{\phantom{11}}&0_{\phantom{12}}\\
0_{\phantom{1}}&0_{\phantom{11}}&0_{\phantom{12}}\\
0_{\phantom{1}}&0_{\phantom{12}}&0_{\phantom{12}}
\end{matrix}
\\
\begin{matrix}
0_{\phantom{1}}&0_{\phantom{11}}&0_{\phantom{12}}\\
0_{\phantom{1}}&0_{\phantom{11}}&0_{\phantom{12}}\\
0_{\phantom{1}}&0_{\phantom{12}}&0_{\phantom{12}}
\end{matrix}
&
\,\,\,\,
\begin{pmatrix}
a&b_1h& b_2h\\
c_1&d_{11}h&d _{12}h\\
c_2&d_{21}h& d_{22}h
\end{pmatrix}^{t-1}
\end{pmatrix}
\begin{pmatrix}
u^+\\
x_1^+\\
x_2^+\\
u^-\\
x_1^-\\
x_2^-\\
\end{pmatrix}
.
$$
or 
$$
\begin{pmatrix}
v^+\\
y_1^+\\
y_2^+\\
v^-\\
y_1^-\\
y_2^-\\
\end{pmatrix}
=
\begin{pmatrix}
\begin{matrix}
a&b_1& b_2\\
c_1&d_{11}&d _{12}\\
c_2&d_{21}& d_{22}
\end{matrix}
&
\begin{matrix}
0_{\phantom{1}}&0_{\phantom{11}}&0_{\phantom{12}}\\
0_{\phantom{1}}&0_{\phantom{11}}&0_{\phantom{12}}\\
0_{\phantom{1}}&0_{\phantom{12}}&0_{\phantom{12}}
\end{matrix}
\\
\begin{matrix}
0_{\phantom{1}}&0_{\phantom{11}}&0_{\phantom{12}}\\
0_{\phantom{1}}&0_{\phantom{11}}&0_{\phantom{12}}\\
0_{\phantom{1}}&0_{\phantom{12}}&0_{\phantom{12}}
\end{matrix}
&
\,\,\,\,
\begin{pmatrix}
a&b_1& b_2\\
c_1&d_{11}&d _{12}\\
c_2&d_{21}& d_{22}
\end{pmatrix}^{t-1}
\end{pmatrix}
\begin{pmatrix}
u^+\\
hx_1^+\\
hx_2^+\\
u^-\\
hx_1^-\\
hx_2^-\\
\end{pmatrix}
.
$$
We introduce new variables
$\wt x^\pm_1=h x^\pm_1$, $\wt x^\pm_2=h x^\pm_2$ and come to the equation
for $\chi_g$. Notice that modules $Q\otimes_{\Z_p} \Z_p^m$ are invariant with respect
to $\OO(m,\Z_p)$.

\sm

{\bf\punct Proof of Proposition \ref{pr:another-quotient}.} 
Proof is the same, we only take $h\in\OO(m,\Q_p)$.
If $Q\subset \cH$  is a subspace, then
 $Q\otimes \ell_m=Q \otimes \Q_p^m$ is a subspace, it is 
 $\OO(m,\Q_p)$-invariant.

\sm


\sm

{\bf\punct  Reformulation of definition.}
The equation (\ref{eq:main-1}) determines a linear subspace
in
$$
\Bigl(\cV\oplus(\cH\otimes \ell_m)\Bigr) \oplus \Bigl(\cV\oplus(\cH\otimes \ell_m)\Bigr)
.
$$
We regard it as a linear relation
$$
\xi:
\Bigl((\cH\otimes \ell_m)\oplus (\cH\otimes \ell_m)\Bigr)
\tto
\bigl(\cV\oplus \cV\bigr)
.$$ 
Then $\chi_\frg$ is the image of the submodule
$$
\eta_{Q,T}=(Q\otimes_{\Z_p} \Z_p^m) \oplus (T\otimes_{\Z_p} \Z_p^m) 
$$
under $\xi$.

\sm


{\bf\punct Immediate corollaries.}
The relation $\xi$ is a morphism of the category $\ov\Naz$. A module
$\eta_{Q,T}$ is self-dual.  By Theorem \ref{th:nazarov-0}
the module $\xi\, \eta_{Q,T}$ is  self-dual.
Theorem \ref{th:self-duality} is proved.

The same argument implies Theorem \ref{th:absolute}.a and Proposition 
\ref{pr:almost-self-dual}.

\sm

Also Lemma \ref{l:choice} became obvious.

\sm


{\bf\punct Continuity (Theorem \ref{th:continuity}).} 
We refer to Theorem \ref{th:continuity0}.

\sm


{\bf\punct Products. Proof of Theorem \ref{th:multiplication}.}
To shorten notation, set $k=2$. Let
{\small
$$
g=
\begin{pmatrix}
a&b_1& b_2\\
c_1&d_{11}&d _{12}\\
c_2&d_{21}& d_{22}
\end{pmatrix}\in\GL(\alpha+2l,\Q_p), \qquad 
h=
\begin{pmatrix}
a'&b'_1& b'_2\\
c'_1&d'_{11}&d'_{12}\\
c'_2&d'_{21}& d'_{22}
\end{pmatrix}
\in\GL(\alpha+2m,\Q_p)
.
$$}
Let $v\oplus w\in\chi_\frg(Q,T)$, $u\oplus v\in \chi_{\frh}(Q,T)$.
Then there are $x\in Q\otimes_{\Z_p}\Z_p^m$, $y\in T\otimes_{\Z_p}\Z_p^m$
such that
{\small
\begin{equation}
\begin{pmatrix}
v^+\\
y_1^+\\
y_2^+\\
v^-\\
y_1^-\\
y_2^-\\
\end{pmatrix}
=
\begin{pmatrix}
\begin{matrix}
a'&b'_1& b'_2\\
c'_1&d'_{11}&d' _{12}\\
c'_2&d'_{21}& d'_{22}
\end{matrix}
&
\begin{matrix}
0_{\phantom{1}}&0_{\phantom{11}}&0_{\phantom{12}}\\
0_{\phantom{1}}&0_{\phantom{11}}&0_{\phantom{12}}\\
0_{\phantom{1}}&0_{\phantom{12}}&0_{\phantom{12}}
\end{matrix}
\\
\begin{matrix}
0_{\phantom{1}}&0_{\phantom{11}}&0_{\phantom{12}}\\
0_{\phantom{1}}&0_{\phantom{11}}&0_{\phantom{12}}\\
0_{\phantom{1}}&0_{\phantom{12}}&0_{\phantom{12}}
\end{matrix}
&
\,\,\,\,
\begin{pmatrix}
a'&b'_1& b'_2\\
c'_1&d'_{11}&d' _{12}\\
c'_2&d'_{21}& d'_{22}
\end{pmatrix}^{t-1}
\end{pmatrix}
\begin{pmatrix}
u^+\\
x_1^+\\
x_2^+\\
u^-\\
x_1^-\\
x_2^-\\
\end{pmatrix}
\label{eq:h-long}
.
\end{equation}}
Also there are $X\in Q\otimes_{\Z_p}\Z_p^l$, $Y\in T\otimes_{\Z_p}\Z_p^l$ 
such that
{\small
\begin{equation}
\begin{pmatrix}
w^+\\
Y_1^+\\
Y_2^+\\
w^-\\
Y_1^-\\
Y_2^-\\
\end{pmatrix}
=
\begin{pmatrix}
\begin{matrix}
a&b_1& b_2\\
c_1&d_{11}&d _{12}\\
c_2&d_{21}& d_{22}
\end{matrix}
&
\begin{matrix}
0_{\phantom{1}}&0_{\phantom{11}}&0_{\phantom{12}}\\
0_{\phantom{1}}&0_{\phantom{11}}&0_{\phantom{12}}\\
0_{\phantom{1}}&0_{\phantom{12}}&0_{\phantom{12}}
\end{matrix}
\\
\begin{matrix}
0_{\phantom{1}}&0_{\phantom{11}}&0_{\phantom{12}}\\
0_{\phantom{1}}&0_{\phantom{11}}&0_{\phantom{12}}\\
0_{\phantom{1}}&0_{\phantom{12}}&0_{\phantom{12}}
\end{matrix}
&
\,\,\,\,
\begin{pmatrix}
a&b_1& b_2\\
c_1&d_{11}&d _{12}\\
c_2&d_{21}& d_{22}
\end{pmatrix}^{t-1}
\end{pmatrix}
\begin{pmatrix}
v^+\\
X_1^+\\
X_2^+\\
v^-\\
X_1^-\\
X_2^-\\
\end{pmatrix}
\label{eq:g-long}
.
\end{equation}}
We write (\ref{eq:g-long}) as
{\footnotesize
$$
\begin{pmatrix}
w^+\\
Y_1^+\\
y^+_1\\
Y_2^+\\
y^+_2\\
w^-\\
Y_1^-\\
y^-_1\\
Y_2^-\\
y^-_2\\
\end{pmatrix}
=
\begin{pmatrix}
\begin{matrix}
a&b_1&0& b_2&0\\
c_1&d_{11}&0&d _{12}&0\\
0&0&1&0&0\\
c_2&d_{21}&0& d_{22}&0\\
0&0&0&0&1
\end{matrix}
&
\begin{matrix}
0_{\phantom{1}}&0_{\phantom{11}}&0_{\phantom{12}}&0_{\phantom{12}}&0\\
0_{\phantom{1}}&0_{\phantom{11}}&0_{\phantom{12}}&0_{\phantom{12}}&0\\
0_{\phantom{1}}&0_{\phantom{12}}&0_{\phantom{12}}&0_{\phantom{12}}&0\\
0_{\phantom{1}}&0_{\phantom{12}}&0_{\phantom{12}}&0_{\phantom{12}}&0\\
0_{\phantom{1}}&0_{\phantom{12}}&0_{\phantom{12}}&0_{\phantom{12}}&0
\end{matrix}
\\
\,\,\,\,\,\,\,
\begin{matrix}
0_{\phantom{1}}&0_{\phantom{1}}&0_{\phantom{1}}&0_{\phantom{1}}&0_{\phantom{1}}\\
0_{\phantom{1}}&0_{\phantom{1}}&0_{\phantom{1}}&0_{\phantom{1}}&0_{\phantom{1}}\\
0_{\phantom{1}}&0_{\phantom{1}}&0_{\phantom{1}}&0_{\phantom{1}}&0_{\phantom{1}}\\
0_{\phantom{1}}&0_{\phantom{1}}&0_{\phantom{1}}&0_{\phantom{1}}&0_{\phantom{1}}\\
0_{\phantom{1}}&0_{\phantom{1}}&0_{\phantom{1}}&0_{\phantom{1}}&0_{\phantom{1}}
\end{matrix}
&
\,\,\,\,\,
\begin{pmatrix}
a&b_1&0& b_2&0\\
c_1&d_{11}&0&d _{12}&0\\
0&0&1&0&0\\
c_2&d_{21}&0& d_{22}&0\\
0&0&0&0&1\\
\end{pmatrix}^{t-1}
\end{pmatrix}
\begin{pmatrix}
v^+\\
X_1^+\\
y_1^+\\
X_2^+\\
y_2^+\\
v^-\\
X_1^-\\
y_1^+\\
X_2^-\\
y_2^-
\end{pmatrix}
.
$$}
Applying
(\ref{eq:h-long}) we come to
{\footnotesize
\begin{multline*}
\begin{pmatrix}
w^+\\
Y_1^+\\
y^+_1\\
Y_2^+\\
y^+_2\\
w^-\\
Y_1^-\\
y^-_1\\
Y_2^-\\
y^-_2\\
\end{pmatrix}
=
\begin{pmatrix}
\begin{matrix}
a&b_1&0& b_2&0\\
c_1&d_{11}&0&d _{12}&0\\
0&0&1&0&0\\
c_2&d_{21}&0& d_{22}&0\\
0&0&0&0&1
\end{matrix}
&
\begin{matrix}
0_{\phantom{1}}&0_{\phantom{11}}&0_{\phantom{12}}&0_{\phantom{12}}&0\\
0_{\phantom{1}}&0_{\phantom{11}}&0_{\phantom{12}}&0_{\phantom{12}}&0\\
0_{\phantom{1}}&0_{\phantom{12}}&0_{\phantom{12}}&0_{\phantom{12}}&0\\
0_{\phantom{1}}&0_{\phantom{12}}&0_{\phantom{12}}&0_{\phantom{12}}&0\\
0_{\phantom{1}}&0_{\phantom{12}}&0_{\phantom{12}}&0_{\phantom{12}}&0
\end{matrix}
\\
\,\,\,\,\,\,\,
\begin{matrix}
0_{\phantom{1}}&0_{\phantom{1}}&0_{\phantom{1}}&0_{\phantom{1}}&0_{\phantom{1}}\\
0_{\phantom{1}}&0_{\phantom{1}}&0_{\phantom{1}}&0_{\phantom{1}}&0_{\phantom{1}}\\
0_{\phantom{1}}&0_{\phantom{1}}&0_{\phantom{1}}&0_{\phantom{1}}&0_{\phantom{1}}\\
0_{\phantom{1}}&0_{\phantom{1}}&0_{\phantom{1}}&0_{\phantom{1}}&0_{\phantom{1}}\\
0_{\phantom{1}}&0_{\phantom{1}}&0_{\phantom{1}}&0_{\phantom{1}}&0_{\phantom{1}}
\end{matrix}
&
\,\,\,\,\,
\begin{pmatrix}
a&b_1&0& b_2&0\\
c_1&d_{11}&0&d _{12}&0\\
0&0&1&0&0\\
c_2&d_{21}&0& d_{22}&0\\
0&0&0&0&1
\end{pmatrix}^{t-1}
\end{pmatrix}
\times\\
\times
\begin{pmatrix}
\begin{matrix}
a'&0&b'_1&0& b'_2\\
0&1&0&0&0\\
c'_1&0&d'_{11}&0&d'_{12}\\
0&0&0&1&0\\
c'_2&0&d'_{21}&0& d'_{22}
\end{matrix}
&
\begin{matrix}
0_{\phantom{1}}&0_{\phantom{11}}&0_{\phantom{12}}&0_{\phantom{12}}&0\\
0_{\phantom{1}}&0_{\phantom{11}}&0_{\phantom{12}}&0_{\phantom{12}}&0\\
0_{\phantom{1}}&0_{\phantom{12}}&0_{\phantom{12}}&0_{\phantom{12}}&0\\
0_{\phantom{1}}&0_{\phantom{12}}&0_{\phantom{12}}&0_{\phantom{12}}&0\\
0_{\phantom{1}}&0_{\phantom{12}}&0_{\phantom{12}}&0_{\phantom{12}}&0
\end{matrix}
\\
\,\,\,\,\,\,\,
\begin{matrix}
0_{\phantom{1}}&0_{\phantom{1}}&0_{\phantom{1}}&0_{\phantom{1}}&0_{\phantom{1}}\\
0_{\phantom{1}}&0_{\phantom{1}}&0_{\phantom{1}}&0_{\phantom{1}}&0_{\phantom{1}}\\
0_{\phantom{1}}&0_{\phantom{1}}&0_{\phantom{1}}&0_{\phantom{1}}&0_{\phantom{1}}\\
0_{\phantom{1}}&0_{\phantom{1}}&0_{\phantom{1}}&0_{\phantom{1}}&0_{\phantom{1}}\\
0_{\phantom{1}}&0_{\phantom{1}}&0_{\phantom{1}}&0_{\phantom{1}}&0_{\phantom{1}}
\end{matrix}
&
\,\,\,\,\,
\begin{pmatrix}
a'&0&b'_1&0& b'_2\\
0&1&0&0&0\\
c'_1&0&d'_{11}&0&d'_{12}\\
0&0&0&1&0\\
c'_2&0&d'_{21}&0& d'_{22}
\end{pmatrix}^{t-1}
\end{pmatrix}
\begin{pmatrix}
u^+\\
X_1^+\\
x_1^+\\
X_2^+\\
x_2^+\\
u^-\\
X_1^-\\
x_1^+\\
X_2^-\\
x_2^-
\end{pmatrix}
\end{multline*}
}
Now 
$$
X\oplus x\,\in\, Q\,\otimes (\Z_p^l\oplus \Z_p^m),\qquad  
Y\oplus y\,\in\, T\otimes (\Z_p^l\oplus \Z_p^m)
,$$
and we get  $u\oplus w\in \chi_{\frg\star\frh}(Q,T)$.
Thus,
$$
\chi_{\frg\star\frh}(Q,T)\supset \chi_{\frg}(Q,T)\chi_{\frh}(Q,T)
.$$
But both sides are self-dual, therefore they coincide. 

\sm


{\bf\punct Morphisms of graphs (Theorem \ref{th:morphisms-graphs}).}
Consider the map 
$$\LMod(\cH)\times \LMod(\cH) \to \LMod(\cH\otimes \ell_m)\times \LMod(\cH\otimes \ell_m) $$
given by $(Q,T)\mapsto (Q\otimes_{\Z_p}\Z_p^m,T\otimes_{\Z_p}\Z_p^m)$.

\begin{lemma}
This map is a morphism of graphs 
$$
\Delta(\cH\Join \cH)\to \Delta\bigl((\cH\otimes \ell_m)\Join (\cH\otimes \ell_m)\bigr)
.$$
\end{lemma}

This statement is obvious.

\sm

Next, we have an embedding of complexes
$$
\Bd\bigl((\cH\otimes \ell_m)\Join (\cH\otimes \ell_m)\bigr)
\to \Bd\bigl((\cH\otimes \ell_m)\oplus (\cH\otimes \ell_m)\bigr)
.
$$

On the other hand,
the linear relation $\xi$ is a morphism of the category $\Naz$. Therefore it
 induces a morphism of graphs
$\Delta\bigl((\cH\otimes \ell_m)\oplus  (\cH\otimes \ell_m)\bigr) \to \Delta(\cV\oplus \cV)$, see
\cite{Ner-gauss}, Proposition 10.7.6.


\sm

{\bf\punct Proof of Proposition \ref{pr:Lambda}.}
We have 
$$
\indef \xi=\Lambda(\frg).
$$
Therefore 
$\Lambda(\frg)\subset\xi\,\eta_{Q,T}\subset \Lambda(\frg)^\bot$. 
This is the statement a) of Proposition \ref{pr:Lambda}.

Also, if $R$ is a relation $\cV\tto W$, $Y\subset \cV$ is a lattice, then
$(RY)_\downarrow=(\indef R)_\downarrow$. This implies b).


\sm

{\bf\punct Values on the distinguished boundary.%
\label{ss:determinant}} Now let $Q$, $T$ be Lagrangian subspaces 
in $\cH$.

\sm

{\sc Proof of Proposition \ref{pr:exotic}}.
Decompose $\cH=\cH_+\oplus \cH_-=\Q_p^\alpha\oplus \Q_p^\alpha$.
A Lagrangian subspace $Q\subset \cH$ of general position is a graph of an
operator $\cH^+\to \cH^-$, and matrix of this operator is symmetric 
(see, e.g., \cite{Ner-gauss}, Theorem 3.1.4). 
To shorten  notation, set $k=2$.
The equation (\ref{eq:main-1}) can be written in the form
\begin{equation}
\begin{pmatrix}
v^+\\
y^+_1\\
y^+_2\\
u^-\\
t_{11}x^+_1+t_{12} x^+_2\\
t_{12}x^+_1+t_{22} x^+_2
\end{pmatrix}
\begin{pmatrix}
a&b_1&b_2& 0&0&0\\
c_1&d_{11}&d_{12}&0&0&0\\
c_2&d_{21}&d_{22}&0&0&0\\
0&0&0&a^t&c^t_1&c^t_2\\
0&0&0&b^t_1&d^t_{11}&d^t_{21}\\
0&0&0&b^t_2&d^t_{12}&d^t_{22}\\
\end{pmatrix}
\begin{pmatrix}
u^+\\
x^+_1\\
x^+_2\\
v^-\\
q_{11}y^+_1+q_{12} y^+_2\\
q_{12}y^+_1+q_{22} y^+_2
\end{pmatrix}
\label{eq:main-2}
,
\end{equation}
We denote
$$
\kappa:=\begin{pmatrix} q_{11}&q_{12}\\q_{12}&q_{22}\end{pmatrix},
\qquad
\tau:=\begin{pmatrix}t_{11}&t_{12}\\t_{12}&t_{22}\end{pmatrix}
$$
and write (\ref{eq:main-2}) as
\begin{align}
&v^+=a u^+ +b x^+
\label{eq:al-1}
\\
&y^+=c u^+ + d x^+
\label{eq:al-2}\\
&u^-=a^t v^-+c^t\kappa y^+
\label{eq:al-3}\\
&\tau x^+= b^t v^-+d^t \kappa y^+
\label{eq:al-4}
.\end{align}
We regard lines (\ref{eq:al-2}),(\ref{eq:al-4})
as a system of equations for $x^+$, $y^+$. The matrix
of the system is
$$
\Omega(\kappa,\tau)=
\begin{pmatrix}
-d&1\\
\tau&-d^t \kappa
\end{pmatrix}
.$$

Evidently, the polynomial $\det\Omega(\kappa,\tau)$ is not zero. Indeed,
fix $\kappa$ and take $\tau=p^{-N}\cdot 1$. If $N$ is sufficiently large,
then the determinant is $\ne 0$. Thus, outside the  hypersurface
$$
\det \Omega(\kappa,\tau)=0
$$
we can express $x^+$ and $y^+$ as functions of $u^+$, $v^-$.
After substitution of $x^+$, $y^+$  to (\ref{eq:al-1}),(\ref{eq:al-3}),
we get a dependence of $u^-$, $v^+$ in $u^+$, $v^-$.
\hfill$\square$

\sm

This also proves Theorem \ref{th:absolute}.b (rationality of characteristic function).

\sm

{\sc Proof Theorem \ref{th:absolute}.c.}
Denote
$$
\begin{pmatrix}
a&b\\c&d
\end{pmatrix}^{-1}=
\begin{pmatrix}
A&B\\C&D
\end{pmatrix}
$$
and write the equation (\ref{eq:main-1}) in the form
$$
\begin{pmatrix}
v^+\\
y^+_1\\
y^+_2\\
v^-\\
q_{11}y^+_1+q_{12} y^+_2\\
q_{12}y^+_1+q_{22} y^+_2
\end{pmatrix}
\begin{pmatrix}
a&b_1&b_2& 0&0&0\\
c_1&d_{11}&d_{12}&0&0&0\\
c_2&d_{21}&d_{22}&0&0&0\\
0&0&0&A^t&C^t_1&C^t_2\\
0&0&0&B^t_1&D^t_{11}&D^t_{21}\\
0&0&0&B^t_2&D^t_{12}&D^t_{22}\\
\end{pmatrix}
\begin{pmatrix}
u^+\\
x^+_1\\
x^+_2\\
u^-\\
t_{11}x^+_1+t_{12} x^+_2\\
t_{12}x^+_1+t_{22} x^+_2
\end{pmatrix}
,$$
or
\begin{align}
v^+= a u^+ + b x^+
\label{eq:1}
\\
y^+= c u^+  +d x^+
\label{eq:2}
\\
v^-= A^t u^- + C^t \tau x^+
\label{eq:3}
\\
 y_+= B^tu^- + D^t \tau x^+
\label{eq:4}
.\end{align}

 We consider lines (\ref{eq:2}), (\ref{eq:4})
 as equations for $y^+$, $x^+$. The matrix of the system is
$$
\Xi(\kappa,\tau)=
\begin{pmatrix}
1&-d\\
\kappa&-D^t \tau
\end{pmatrix}
.$$
Its determinant equals
$$\det\Xi(\kappa,\tau)=\det(-D^t \tau+\kappa d).$$
If it is nonzero, we get a linear operator $u\mapsto v$.
 We come to the following statement:

\begin{proposition}
If there exists a pair of symmetric matrices $\kappa$, $\tau$
such that $\det(-D^t \tau+\kappa d)\ne 0$, then $\chi_\frg(Q,T)\in\Sp(\cV,\Q_p)$
a.s. on $\LGr(\cH)\times\LGr(\cH)$.
\end{proposition}


{\bf\punct Involution. Proof of Proposition \ref{pr:involution}.}
 We write the defining relation for $\chi_{g^{-1}}$,
 $$
 \begin{pmatrix}
 v^+\\
 y^+\\
 v^-\\
 y^-
 \end{pmatrix}
 =
 \begin{pmatrix}
 \begin{pmatrix} a&b\\c&d\end{pmatrix}^{-1}&
\phantom{ \biggl(} \begin{matrix} 0&0\\0&0\end{matrix}\phantom{\biggr)^{t}}\\
 \phantom{ \biggl(}  \begin{matrix} 0&0\\0&0\end{matrix} \phantom{\biggr)^{t}}&
  \begin{pmatrix} a&b\\c&d\end{pmatrix}^t 
 \end{pmatrix}
  \begin{pmatrix}
 u^+\\
 x^+\\
 u^-\\
 x^-
 \end{pmatrix}
, $$
 represent this in the form
 $$
   \begin{pmatrix}
 u^+\\
 x^+\\
 u^-\\
 x^-
 \end{pmatrix}
=
 \begin{pmatrix}
 \begin{matrix} a&b\\c&d\end{matrix}&
\phantom{ \biggl(} \begin{matrix} 0&0\\0&0\end{matrix}\phantom{\biggr)^{t-1}}\\
  \begin{matrix} 0&0\\0&0\end{matrix} &
  \begin{pmatrix} a&b\\c&d\end{pmatrix}^{t-1}
 \end{pmatrix}
 \begin{pmatrix}
 v^+\\
 y^+\\
 v^-\\
 y^-
 \end{pmatrix}
 $$
and come to desired statement.


\sm

{\bf\punct  Proof of Theorem \ref{th:symmetry}.}
We write (\ref{eq:main-1}) as
{\footnotesize
\begin{multline*} 
    \begin{pmatrix}
 u^+\\
 x^+\\
 u^-\\
 x^-
 \end{pmatrix}
=
 \begin{pmatrix}
 \lambda^{-1}&&&\\
 &\lambda^{-1}&&\\
 &&\lambda&\\
 &&&\lambda
 \end{pmatrix}
 \begin{pmatrix}
 \begin{matrix} a&b\\c&d\end{matrix}&
\phantom{ \biggl(} \begin{matrix} 0&0\\0&0\end{matrix}\phantom{\biggr)^{t-1}}\\
  \begin{matrix} 0&0\\0&0\end{matrix} &
  \begin{pmatrix} a&b\\c&d\end{pmatrix}^{t-1}
 \end{pmatrix}
 \begin{pmatrix}
 \lambda&&&\\
 &\lambda&&\\
 &&\lambda^{-1}&\\
 &&&\lambda^{-1}
 \end{pmatrix}
 \begin{pmatrix}
 v^+\\
 y^+\\
 v^-\\
 y^-
 \end{pmatrix}
 \end{multline*}}
or
$$ 
    \begin{pmatrix}
\lambda u^+\\
\lambda x^+\\
\lambda^{-1} u^-\\
\lambda^{-1} x^-
 \end{pmatrix}
=
 \begin{pmatrix}
 \begin{matrix} a&b\\c&d\end{matrix}&
\phantom{ \biggl(} \begin{matrix} 0&0\\0&0\end{matrix}\phantom{\biggr)^{t-1}}\\
  \begin{matrix} 0&0\\0&0\end{matrix} &
  \begin{pmatrix} a&b\\c&d\end{pmatrix}^{t-1}
 \end{pmatrix}
 \begin{pmatrix}
\lambda v^+\\
\lambda y^+\\
\lambda^{-1} v^-\\
\lambda^{-1} y^-
 \end{pmatrix}
$$

{\bf\punct Another reformulation of the definition of characteristic functions.%
\label{ss:another}}
Consider the space $W=\cV\oplus (\cH\otimes \ell_m)$. For any self-dual submodule
$Q\subset \cH$,
consider the linear relation $\Lambda:\cV\tto W$ defined by
$$
\Lambda_Q=1_\cV\oplus (Q\otimes \Z_p^m)\,\subset\, (\cV\oplus \cV)\oplus (Q\otimes \ell_m)
.
$$
Then $\chi_\frg$ is a product of linear relations 
$$
\chi_\frg(Q,T)=(\Lambda_T)^\square \begin{pmatrix} g&0\\0& g^{t-1}\end{pmatrix} \Lambda_Q
.
$$




\section{Multiplicativity theorem}

\COUNTERS

Theorem \ref{th:multiplicativity} (multiplicativity theorem) formulated above
is a representative of wide class of theorems, their proofs are standard,
below we refer to  proofs \cite{Ner-book}, Chapter VIII.

\sm


{\bf\punct Corners of orthogonal matrices.}

\begin{lemma}
Let $A$ be a $m\times m$ matrix with  elements $\in\Z_p$.
Then there exists $N$ and a matrix $\begin{pmatrix} A&B\\C&D\end{pmatrix}\in \OO(m+N,\Z_p)$.
\label{l:corner}
\end{lemma}

{\sc Proof.} Denote by $\bfB_m$ the set of all possible $m\times m$ left upper corners
of matrices $g\in \OO(\infty,\Z_p)$. 

1) The set $\bfB_m$ is closed with respect to matrix products. Indeed, let 
$$
\begin{pmatrix} A&B\\C&D\end{pmatrix} \in \OO(m+N,\Z_p),\qquad
 \begin{pmatrix} A'&B'\\C'&D'\end{pmatrix} \in \OO(m+N',\Z_p)
. $$
Then
$$
\begin{pmatrix} A&B&0\\C&D&0\\0&0&1 \end{pmatrix}
 \begin{pmatrix} A'&0&B'\\0&1&0\\C'&0&D'\end{pmatrix}=
 \begin{pmatrix}
 AA'&\dots&\dots\\
\dots &\dots&\dots\\
\dots &\dots&\dots
 \end{pmatrix}\in \OO(m+N+N',\Z_p).
$$

2) If $A\in \bfB_m$, $A'\in\bfB_n$, then 
$\begin{pmatrix}A&0\\0&A'\end{pmatrix}\in \bfB_{m+n}$.

\sm

3) It is more-or-less clear that for any  $z\in\Z_p$ we have
$$
\begin{pmatrix} z\end{pmatrix}\in \bfB_1,\qquad
\begin{pmatrix}1& z\\0&1\end{pmatrix}, \begin{pmatrix}1& 0\\z&1\end{pmatrix}\in\bfB_2.
$$

\sm

4) $\bfB_m$ contains matrices of permutations.

\sm

Now we can produce any matrix with integer  elements.
\hfill $\square$


\sm

{\bf\punct Admissible representations.}
Denote by $\bfK_m$ the subgroup in $\bfK$ consisting of matrices of the form
$\begin{pmatrix} 1_m&0\\0&*\end{pmatrix}$. 

Let $\tau$ be  a unitary representation of $\bfK$ in a Hilbert space $H$.
Denote by $H(m)$ the subspace of $\bfK_m$-fixed vectors. Denote by $P(m)$
the operator of orthogonal projection to $H(m)$. We say, that
$\tau$ is {\it admissible} if $\cup_m H(m)$ is dense in $H$.

We say, that a representation of $\bfG$ is {\it $\bfK$-admissible} if
its restriction to $\bfK$ is admissible.



\sm


{\bf\punct Continuation of representations.}
Denote by $\bfB_\infty$ the semigroup of all infinite matrices $A$ such that:

\sm

a) $a_{ij}\in \Z_p$;

\sm

b) for each $i$ the sequence $a_{ij}$ tends to $0$ as $j\to\infty$;
for each $j$ the sequence $a_{ij}$ tends to $0$ as $i\to\infty$.

\sm
  
We say that a sequence of matrices $A^{(j)}\in \bfB_\infty$ weakly converges to $A$ 
if we have convergence of each matrix element,   $a_{kl}^{(j)}\to a_{kl}$.

\sm

Denote by $\bfO(\infty,\Z_p)$ the group of all orthogonal matrices
$\in \bfB_\infty$. 

\begin{lemma}
The group $\OO(\infty,\Z_p)$ is dense in $\bfO(\infty,\Z_p)$ and in $ \bfB_\infty$.
\end{lemma}

{\sc Proof.} Let $S\in \bfB_\infty$. Consider its left upper corner of  size
$m\times m$. Consider $g_m\in \O(\infty,\Z_p)$ having the same left upper corner.
Then $g_m$ weakly converges to $S$,
\hfill $\square$

\begin{theorem}
{\rm a)} Let $\tau$ be a unitary representation of $\bfK=\OO(\infty,\Z_p)$.
 The following conditions are equivalent:

\sm

--- $\tau$ is admissible;

\sm

--- $\tau$ admits a weakly continuous extension to the group  $\bfO(\infty,\Z_p)$;

\sm

--- $\tau$ admits a weakly continuous extension to a representation $\wt\tau$
of the semigroup $\bfB_\infty$ such that
$\wt\tau(A^t)=\wt\tau(A)^*$, $\|\wt\tau(A)\|\le 1$ for all $A$.

\sm

{\rm b)} For an admissible representation $\tau$,
 $$P(m)=\wt\tau \begin{pmatrix}1_m&0\\0&0\end{pmatrix}.$$
 \end{theorem}
 
This is a statement in the spirit of \cite{Olsh-kir}. We omit a proof, since 
it is a one-to-one repetition of proof of \cite{Ner-book}, Theorem VIII.1.4
about symmetric groups (admissibility implies semigroup continuation),
the only new detail is Lemma \ref{l:corner}). Admissibility follows from
continuity by \cite{Ner-book}, Proposition VIII.1.3. 
 
 \begin{corollary}
 \label{cor:Theta}
 Denote
 $$
 \Theta_N^{(m)}= 
 \begin{pmatrix}
1_m&0&0&0\\
0&0&1_N&0\\
0&1_N&0&0\\
0&0&0&1_\infty
\end{pmatrix} 
.
$$
The projector $P(m)$ is a 
 weak limit of the sequence
 \begin{equation}
 P(m)=
 \lim_{N\to\infty}\tau(\Theta_N^{(m)})
   \label{eq:lim}
   .
   \end{equation}
\end{corollary}

{\sc Proof.} The sequence $\Theta_N^{(m)}\in \OO(\infty,\Z_p)$ weakly converges
to the matrix $\begin{pmatrix}1_m&0\\0&0\end{pmatrix}\in\bfB_\infty$. We refer to the statement
b) of the theorem.
\hfill $\square$


\sm

{\bf\punct Proof of Theorem \ref{th:multiplicativity}.}
We keep the notation of Subsection \ref{ss:multiplicativity}.
 Let $v\in H^\bfK$, $g\in G_j=\GL(\alpha+km,\Q_p)$,
let $q\in \bfK_j$. Then
$$
\rho(q)\rho(g) v= \rho(g)\rho(q)h= \rho(g)h
,$$
i.e., $v\in H(j)$. Thus the subspace $\cup_j H(j)$ is $\bfG$-invariant.
Its closure is an admissible representation of $\bfG$. 
In $\bigl(\cup_j H(j)\bigr)^\bot$ Theorem \ref{th:multiplicativity}
holds by a trivial reason (the space of fixed vectors $\bfK$ is zero).

Thus, without loss of generality we can assume that $\rho$ is admissible.

Now let $g$, $h\in \bfG$, let $\frg$, $\frh\in\bfK\setminus \bfG/\bfK$ be the corresponding double cosets. 
Let $P=P(0)$ be the projector to $\bfK$-fixed vectors.
  Applying Corollary \ref{cor:Theta},
we obtain
$$
\ov\rho(\frg)\ov\rho(h)=
P\rho(g)P\rho(h)=
\lim_{N\to\infty}
P\rho(g)\rho(\frI(\Theta_N^{(0)}))\rho(h)
=\lim_{N\to\infty}
P\rho(g\frI(\Theta_N) h)
,
$$
here $\frJ:\bfK\to\bfG$ is the embedding (\ref{eq:tau}).
By the definition ($\Theta_N^{(0)}$ is $\Theta_N$ from Subsection \ref{ss:multiplicativity}), we get $\ov\rho(\frg\star\frh)$.


\sm

{\bf\punct Variation of construction. Train.}
We can define multiplication of double cosets
$$
\bfK_p\setminus \bfG/\bfK_q\,
\times\,
 \bfK_q\setminus \bfG/\bfK_r\,\to \, 
 \bfK_p\setminus \bfG/\bfK_r
.$$
In the definition of product of double cosets (Subsection  \ref{ss:double-cosets}),
we simply change $\Theta_N$ by $\Theta_N^{(q)}$. An explicit formula of the product
is the same (\ref{eq:coset-product}). Thus we get a category ({\it train
$\cT(\bfG,\bfK)$ of 
the pair $(\bfG,\bfK)$}).

Next, for any unitary representation $\rho$
of the group $\bfG$, a double coset $\frg\in \bfK_p\setminus \bfG/\bfK_q$ determines an
operator
$
\ov\rho(\frg):H(q)\to H(p)
$
by  the formula
$$
\ov\rho(g):=P(q)\rho(g),\qquad g\in\frg
.
$$
For any 
$$
\frg\in\bfK_p\setminus \bfG/\bfK_q\qquad
\frh\in
 \bfK_q\setminus \bfG/\bfK_r
,$$
the following identity holds
$$
\rho(\frg)\rho(\frh)=\rho(\frg\star\frh),
$$
i.e., we get a representation of the category $\cT(\bfG,\bfK)$. Also,
 \begin{equation}
 \rho(\frg^*)=\rho(\frg)^*,\qquad
 \|\rho(\frg)\|\le 1
 .
 \label{eq:||}
 \end{equation}

Also it can be shown that

\begin{theorem}
 This construction is a bijection
between the set of $\bfK$-admissible unitary representations of $\bfG$
and the set of representations of the category $\cT(\bfG,\bfK)$
satisfying
{\rm(\ref{eq:||})}.
\end{theorem}

We omit a proof, since it is the same as in \cite{Ner-faa}.
\hfill $\square$

\sm

Also the construction of characteristic functions and their properties
survive for double cosets $\bfK_p\setminus \bfG/\bfK_q$.


\section{Representations of the group  $\bfG$}

\COUNTERS


\sm

{\bf\punct Existence of representations.} Let
$$
\begin{pmatrix}
a&b_1&\dots& b_k\\
c_1&d_{11}&\dots&d _{1k}\\
\vdots&\vdots & \ddots&\vdots\\
c_k&d_{k1}&\dots& d_{kk}
\end{pmatrix}
\in\GL(\alpha+k\infty,\Q_p)
.
$$
Consider embedding
$\GL(\alpha+k\infty,\Q_p) \to \Sp(2(\alpha+k\infty),\Q_p)$
given   by
$$
\iota:
g\mapsto\begin{pmatrix}
g&0\\0&g^{t-1}
\end{pmatrix}
.
$$

For any 
$$
r=
\begin{pmatrix}
r_{11}&\dots &r_{1\,2n}\\
\vdots &\ddots&\vdots\\
r_{2n\,1}&\dots &r_{2n\,2n}
\end{pmatrix}\in\Sp(2k,\Q_p)
$$
consider the matrix $\sigma(r)=1_{2\alpha}\oplus (r\otimes 1_\infty)$,
$$
\sigma(r):=
\begin{pmatrix}
1_\alpha&0& \dots& 0&0\\
0&r_{11}\cdot 1_\infty& \dots &0& r_{1k}\cdot 1_\infty\\
\vdots&\vdots&\ddots&\vdots&\vdots\\
0&0& \dots& 1_\alpha&0\\
0&r_{11}\cdot 1_\infty& \dots &0& r_{1k}\cdot 1_\infty\\
\end{pmatrix}
$$

This matrix is not contained in $\Sp(2(\alpha+k\infty),\Q_p)$,
because it is not finitary. However, the map
\begin{equation}
q\mapsto \sigma(r^{-1})\, q\,\sigma(r)
\label{eq:automorphism}
\end{equation}
is an outer automorphism  of $\Sp(2(\alpha+k\infty),\Q_p)$. 
Emphasize that this automorphism fixes the subgroup $\bfK=\OO(\infty,\Z_p)$.

We consider the representation $\rho(r)$ of $\GL(\alpha+k\infty,\Q_p)$  given by the formula
$$
\rho_r(g)=\We\bigl(\sigma(r^{-1}) \iota(g)\sigma(r)\bigr)
,$$
where $\We(\cdot)$ is the Weil representation, see Subsection \ref{ss:weil}.

Recall that the Weil representation is projective.

\begin{lemma} 
The representation $\rho_r$ is equivalent to a linear representation, i.e.,
there is a function {\rm(a trivializer)} $\gamma:\bfG\to\C^\times$ such that 
$\gamma(g)\rho_r(g)$ is a linear representation. 
\end{lemma}

{\sc Proof.} First, the restriction of the Weil representation 
of $\Sp(2n,\Q_p)$ to $\GL(n,\Q_p)$ is linear, see (\ref{eq:weil-gl}).
 Therefore, restricting the Weil representation
to each finite-dimensional group $G_j=\GL(\alpha+kj,\Q_p)$ we get a representation
equivalent to a linear representation (for finite-dimensional groups 
the automorphism (\ref{eq:automorphism}) is inner).
Denote by $\gamma_j(g)$ the trivializer for $G_j$.
Ratio $\gamma(g)_j/\gamma(g)_{j+1}$ of two trivializers is a
character  $G_j\to\C^\times$. All characters of $G_j\to\C^\times$
has the form $\phi(\det h)$, where $\phi$ is a character $\Q^\times\to \C^\times$.
Correcting $\gamma_{j+1}(g)\mapsto \gamma_{j+1}(g) \psi(\det g)$, we can assume
that $\gamma_{j+1}(g)=\gamma_j(g)$ on $G_j$.

In this way we choose a trivializer $\gamma$ on the whole group $\bfG$. Restriction of
$\gamma$ to $\OO(\infty,\Z_p)$ must be a character on 
$\OO(\infty,\Z_p)\to\C^\times$. The only non-trivial character is $\det(u)=\pm 1$.
We change the trivializer $\gamma(g)$ to $\det(g)\gamma(g)$.
\hfill $\square$

\begin{lemma}
\label{l:ergodic}
 In the model of Subsection {\rm \ref{ss:weil}}, the 
 subspace $L^2(\cE_{\alpha+k\infty})^\bfK$ of $\bfK$-fixed vectors of $\rho_r$  coincides
with the space of functions of the form
$$
f(z_1,\dots,z_\alpha)I(z_{\alpha+1}) I(z_{\alpha+2})\dots
$$
\end{lemma}

{\sc Proof.}  Without loss of generality, we can set $\alpha=0$.
 We  regard $\cE_{k\infty}$
as the space
of $ \infty \times k$ matrices $Z=\{z_{ij}\}$ with elements in $\Q_p$
(all but a finite number of matrix elements are in $\Z_p$).
The group $\bfK=\OO(\infty,\Z_p)$ acts by left multiplications
$$
\We(u)f(Z)= f(Zu).
$$

We must show that $\prod_{ij} I(z_{ij})$ is a unique $\OO(\infty,\Z_p)$-invariant
function in $L^2(\cE_{k\infty})$. Equivalently, $\Z_p^{k\infty}$ is a 
unique invariant subset of finite positive measure.

The group $\OO(\infty,\Z_p)$ contains the group  $S(\infty)$
of finitely supported permutations of the set $\N$. According
zero-one law (see, e.g., \cite{Shi}, \S4.1), the action of $S(\infty)$ on the set
$\Z_p^{k\infty}\subset \cE_{k\infty}$ is ergodic.
Let $\Omega\subset \cE_{k\infty}$ be an invariant set.
Let $\xi\in \cE_{k\infty}\setminus \Z_p^{k\infty}$.
Assume that the measure of the set $\Omega\cap (\xi+\Z_p^{k\infty})$ is non-zero, say $\nu_0$.
Since $\Omega$ is $S(\infty)$-invariant, for any $\frs\in S(\infty)$,
the set $\Omega\cap (\xi \frs+\Z_p^{k\infty})$ has the same measure 
$\nu_0$. However there is a countable number of disjoint sets 
of the form
$\xi \frs+\Z_p^{k\infty}$, therefore the measure of
$\Omega$ is infinite.
\hfill $\square$

\begin{corollary}
Let $\alpha=0$. Then the representation $\rho_r$ contains a unique
irreducible $\bfK$-spherical representation of $\bfG$.
\end{corollary}

{\sc Proof.} We take the $\bfG$-cyclic span of the unique $\bfK$-fixed vector.
\hfill $\square$

\sm

Next, consider the subgroup $\GL(1,\Q_p)\subset \Sp(2k,\Q_p)$ consisting
of matrices $\begin{pmatrix}\lambda\cdot 1_k&0\\0&\lambda^{-1}\cdot 1_k \end{pmatrix}$,
where $\lambda\in \Q_p^\times$.

\begin{lemma}
If $r$, $r'\in\Sp(2k,\Q_p)$ are contained in the same  double coset
$$\GL(1,\Q_p)\setminus \Sp(2k,\Q_p)/  \Sp(2k,\Z_p),$$
 then $\rho_r\simeq\rho_{r'}$.
\end{lemma}

{\sc Proof.} First, if $q\in \GL(1,\Q_p)$, then the automorphism 
 (\ref{eq:automorphism}) fixes the subgroup $\GL(\alpha+k\infty,\Q_p)$.

Second, if $t\in  \Sp(2k,\Z_p)$, then $\sigma(t)$ is contained
in the group $\mathbf{Sp}$ of automorphisms of the infinite object of the Nazarov category.
Therefore the operator 
$\We(\sigma(t))$ is well-defined, it intertwines  $\rho_r$ and $\rho_{rt}$.
\hfill $\square$

\sm


{\bf\punct Relation of characteristic functions and representations.}
By Lemma \ref{l:ergodic}, we can identify
the space of $\bfK$-fixed vectors of $\rho_r$ and the space of the Weil representation
of $\Sp(2\alpha,\Q_p)$.

\begin{theorem}
\label{th:link} The representation of the semigroup
$\bfK\setminus \bfG/\bfK$ in the space of $\bfK$-fixed vectors of $\rho_r$ 
is given by the formula
$$
\ov\rho_r(\frg)=s\cdot \We\bigl(  \chi_\frg(r\Z_p^{2k}, r\Z_p^{2k})\bigr),\qquad s\in\C^\times
.$$
\end{theorem}

{\sc Proof.} We use the notation and statements  of Subsection \ref{ss:weil}.
Let $\cV$ and $\cH$ be the same as in Section 4.
Let  $Y=\cV_{2k\infty}$, $W=\cV\oplus Y$.
The operator of projection $\cH(\cV\oplus Y)$ to
 $\cH(V\oplus Y)^\bfK\simeq \cH(V)$ is 
$\We(\theta_W^V)$. Therefore
$$
\ov\rho(\frg)=s'\cdot\We(\theta_W^V)\We(\sigma(r^{-1})\iota(g)\sigma(r))\We(\theta_W^V)
$$
as an operator $L^2(\cE_{\alpha+k\infty})^\bfK\to L^2(\cE_{\alpha+k\infty})^\bfK$. The operator
$$\We(\lambda_W^V):L^2(\Q_p^\alpha)\to L^2(\cE_{\alpha+k\infty}) $$
 is an operator of isometric 
embedding, the image is $\cH(V\oplus V_{2k\infty})^\bfK$. Therefore we can write
$\ov\rho(\frg)$ as
\begin{multline}
\ov\rho(\frg)=
s^{''}\cdot\We(\lambda_W^V)^*\We(\theta_W^V)\We(\sigma(r^{-1})\iota(g)\sigma(r))\We(\theta_W^V)\We(\lambda_W^V)
=\\=
s^{'''}\cdot\We(\lambda_W^V)^*\We(\sigma(r^{-1})\iota(g)\sigma(r))\We(\lambda_W^V)=\\=
s^{''''}\cdot\We\Bigl[(\lambda_W^V)^* \sigma(r^{-1})\iota(g)\sigma(r)\lambda_W^V\Bigr]
\label{eq:brackets}
.
\end{multline}
Next, $\sigma(r)\lambda_W^V:V\tto V\oplus Y$ is a direct sum
of $1_V\subset V\oplus V$ and the lattice in $Y$ given by
$$\sigma(r) Y(\O)
=\sigma(r) ( H(\O)\otimes \O^\infty)=
(r H(\O))\otimes \O^\infty)
.$$
We apply Subsection \ref{ss:another} for the expression in square brackets
in (\ref{eq:brackets}).

\sm


{\bf\punct A more general construction.} 
Consider the embedding
$$
\iota_l:\GL(\alpha+k\infty,\Q_p)\to\Sp(2 l\alpha+2lk\infty, \Q_p)
$$
given by
$$
g\mapsto
\begin{pmatrix}
g&\dots&0 &0&\dots &0\\
\vdots&\ddots&\vdots & \vdots&\ddots&\vdots\\
0&\dots&g &0&\dots &0\\
0&\dots&0 &g^{t-1}&\dots &0\\
\vdots&\ddots&\vdots & \vdots&\ddots&\vdots\\
0&\dots&0 &0&\dots &g^{t-1}\\
\end{pmatrix}
.
$$
This is a $2l\times 2l$ block matrix,
 each block of this matrix has size $(\alpha+k\infty)\times (\alpha+k\infty)$.
 
 Next, for a matrix $r\in\Sp(2kl,\Q_p)$ we take
 $$
 \sigma(r):=1_{2\alpha l}\oplus ( r\otimes 1_\infty)
 $$
 and consider the representation of $\GL(\alpha+k\infty,\Q_p)$ given by 
 $$
 \rho_r(g)=\We(\sigma(r)^{-1}\iota_l(g) \sigma(r))
 .$$

Set $\alpha=0$.
As above, each representation $\rho_r$ of $\bfG=\GL(k\infty,\Q_p)$
contains a unique $\bfK$-spherical subrepresentation.

\begin{conjecture}
Any $\bfK$-spherical representation of $\GL(k\infty,\Q_p)$
is a subrepresentation in
$\phi(\det(g))\,\rho_r(g)$, where $\phi=\phi_r:\Q_p^\times\to\C^\times$ is a character.  Representations
$\rho_r$ are parametrized by the set
\begin{equation*}
\bigcup_l
\GL(l,\Q_p)\setminus \Sp(2kl,\Q_p)/ \Sp(2kl,\Z_p).
\end{equation*}
\end{conjecture}



{\tt Math.Dept., University of Vienna,

 OscarMorgensternplatz, 1,
Vienna, Austria

\&

Institute for Theoretical and Experimental Physics,

Bolshaya Cheremushkinskaya, 25, Moscow 117259,
Russia

\&

Mech.Math.Dept., Moscow State University,

Vorob'evy Gory, Moscow

e-mail: neretin(at) mccme.ru

URL:www.mat.univie.ac.at/$\sim$neretin

wwwth.itep.ru/$\sim$neretin
}

\end{document}